\documentclass[largeformat]{interact}

\usepackage[caption=false]{subfig}% Support for small, `sub' figures and tables

\usepackage[numbers,sort&compress]{natbib}% Citation support using natbib.sty
\bibpunct[, ]{[}{]}{,}{n}{,}{,}% Citation support using natbib.sty
\theoremstyle{plain}% Theorem-like structures provided by amsthm.sty

\theoremstyle{definition}

\theoremstyle{remark}

\begin{document}

\articletype{}% Specify the article type or omit as appropriate

\title{A statistical methodology to select covariates in high-dimensional data under dependence. Application to the classification of genetic profiles in oncology.}

\author{
\name{B. Bastien\textsuperscript{c}, T. Boukhobza\textsuperscript{b}, H. Dumond\textsuperscript{b}, A. G\'egout-Petit\textsuperscript{a}, A. Muller-Gueudin\textsuperscript{a}\thanks{A. Muller-Gueudin.
Email: Aurelie.gueudin@univ-lorraine.fr}, and C. Thi\'ebaut\textsuperscript{b}}
\affil{\textsuperscript{a}Universit\'e de Lorraine, CNRS, Inria, IECL, F-54000 Nancy, France; \textsuperscript{b}Universit\'e de Lorraine, CNRS, CRAN, F-54000 Nancy, France;  \textsuperscript{c}Transgene S.A., Boulevard Gonthier d'Andernach, Parc d'innovation, CS80166, F-67405 Illkirch-Graffenstaden Cedex, France} }

\maketitle

\begin{abstract}
We propose a new methodology for selecting and ranking covariates associated with a variable of interest in a context of high-dimensional data under dependence but few observations. The methodology successively intertwines the clustering of  covariates, decorrelation of covariates using Factor Latent Analysis, selection using aggregation of adapted methods and finally ranking.  A   simulation study shows the interest of the decorrelation inside the different clusters of covariates. 
 We first apply our method to transcriptomic data of 37 patients with advanced non-small-cell lung cancer who have received chemotherapy,  to select the transcriptomic covariates that explain the survival outcome of the treatment. Secondly, we apply our method to 79 breast tumor samples to define patient profiles for a new metastatic biomarker and associated gene network in order to personalize the treatments.\end{abstract}


\section{Introduction}
The purpose of personalized medicine is to select appropriate and optimal therapies based on the context of a patient's genetic content or other molecular or cellular analysis. One of the main challenges faced by biologist and mathematician consortium for the construction of explanatory models of multivariable biological processes, is the relatively low amount of experimental data available compared to the huge number of variables. The point is of great concern when the biological question deals with transcriptomic data in order to build gene networks and decipher the role of a rare isoform, for which no specific probe is currently available. In that context, the purpose of this paper is to propose a method  to select the covariates that are linked with the outcome of a given therapy or a biological marker, among a set of more than tens of thousands covariates. For instance, the relevant dataset we will study in this paper are the following: 
\begin{itemize}
\item 51336  transcriptomic data of 37 patients with advanced non-small-cell lung cancer who have received chemotherapy. The survival time being known, the objective is  to select the transcriptomic covariates that explain the survival outcome of the treatment, then to define the profiles of the patients who survive the treatment.
\item 54676 probes in 79 breast tumor samples. A new metastatic biomarker being known, the objective is to define patient profiles for this metastatic biomarker, and associated gene network in order to personalize the treatments.
\end{itemize}
The variable of interest being known (treatment outcome or biological marker), the question is to find its link with the transcriptomic profile of the patients. We propose a methodology, that, firstly selects and ranks the transcriptomic covariates that are the most linked with the outcome treatment, and secondly, that visualises the profiles of the selected transcriptomic covariates, for all the patients of the study.

More generally, the problem to detect association between a variable of interest and many covariates has been tackled by many biologists and statisticians \cite{dudoit2002statistical,aubert2004determination,bar2005comparaisons, le2011sparse, gunther2014novel}.  A common example, coming from biology,
is testing which of $p$ genes' expression levels given in a dataset $\mathbf X$ is linked significantly with a variable $Y$, which we will call the variable of interest. The  variable of interest may  be a binary variable like
an outcome of  treatment or it may be a quantitative variable such as a phenotype or physiological parameter. In the two data studies of this paper, the aim of the biologist is not necessarily to detect exhaustively all the genes involved in his problem but to have a list of the most important of them in order to study their biological functions. For this purpose, it is interesting to rank the genes according to the strength of their link with the variable of interest. 
Although we present biological studies, our goal is to propose a general methodology in a context of high dimensional data (the number $p$ of covariates is in the order of thousands) while the total number $n$ of samples could be small (for instance between $25$ and $100$).

In the context of transcriptomic data, the covariates are high dimensional and correlated. This correlation between covariates, in a high-dimensional context, has to be taken into account in the statistical analysis. Moreover, we are in a context of small sample size ($n\ll p$). Thus, robustness of the statistical analysis has to be quantified.

We cite here some statistical methods that have been developed to select covariates in high-dimensional contexts. The state of the art about the control of false discoveries in multiple testing procedures is very extensive. The famous correction proposed by  \cite{bonferroni1936teoria} to control the Family Wise Error Rate (FWER) has been emulated and we can find a review about these methods in \cite{dudoit2004multiple}. Alternative methods focused on the control of the False Discovery Rate (FDR) (\cite{benjamini1995controlling,benjamini2001control}), or of the local FDR \cite{efron2005local} or the q-value \cite{storey2002direct,storey2003statistical,storey2003positive}. For a  review (in french) of the methods, see \cite{bar2005comparaisons}. Regarding regression in the framework of high dimensional data ($n \ll p$), many methods are available. For exemple, the PLS approach of \cite{tenenhaus2005pls} is a kind of principal component regression. The Lasso regression proposed by \cite{tibshirani1996regression} performs both variable selection and regularization in penalizing the sums of squares by the $L_1$-norm of the coefficients. This method has been derived for many kinds of problems like logistic-regression in the case of binary data  \cite{meier2008group}, or network inference \cite{meinshausen2006high}. Another versatile tool to select covariates  in different non parametric contexts is given by the random forests, with the concept of importance of covariates \cite{genuer2010variable}. 

Another important characteristic of the data that has to be taken into account  is the structure of covariance of the covariates. Most of the multiple testing corrections make the assumption of the independence between the tests. However it is well-known that omics data are correlated by blocks. In the context of multiple testing, covariance between the covariates could bias the uniform repartition of the p-values under the null hypothesis and also inflates the variance of the estimation of the FDR \cite{friguet2009factor,friguet2011estimation}. In \cite{genuer2010variable} it is also shown that despite the robustness of the random forests, importance of covariates calculated by random forests  is perturbed by adding other correlated covariates. One of the ways to deal with dependence is to model it by latent factors; it is a way to reduce the information in supposing that the common information of the $p$ covariates is given by $q \ll p$ latent factors as   \cite{friguet2009factor} and \cite{friguet2011estimation}. More precisely, they  propose a way to correct the data according to a regression link with the variable of interest $Y$ in such a way that covariates are independent conditionally to $Y$ (leading to the independence of the tests). After this correction, they propose a multiple testing procedure based on the  method of \cite{benjamini1995controlling} and \cite{benjamini2001control}.  This method of correction will be called FAMT correction (for Factor Analysis for Multiple Testing) in the sequel.

However, the framework of FAMT is to consider the data $\mathbf X$ as an only one block of correlated covariates and has to be adapted if $\mathbf X$ is structured in several independent clusters of correlated covariates.  As we will see in Section \ref{sec:simu}, the FAMT does not give good results if it is applied directly on the whole set of data $\mathbf X$, without taking into account its decomposition in clusters with strong within correlation. Then, we propose to identify the clusters of correlated covariates before performing FAMT correction on each of the clusters. The clustering of covariates as proposed by  \cite{chavent2011clustofvar} is a good way to arrange covariates into homogeneous clusters, i.e., groups inside of which  covariates are strongly related to each other.   %Bootstrap \cite{efron1992bootstrap} is commonly used in many circumstances particularly to estimate variability when the number of samples is low. It is also used for ranking sometimes trough p-values \cite{mukherjee2003gene,hall2009using}. 

Our purpose in this paper is to propose a method adapted to the selection (and ranking) of correlated quantitative covariates associated with a variable of interest. For this, we propose a methodology that takes into account (1) the structure of correlation by clusters of covariates; %by using a clustering of covariates in order to identify the groups; 
(2) the correlation inside each cluster of correlated covariates.% by applying the analysis in factors proposed by Friguet {\it et al.} in \cite{friguet2009factor, friguet2011estimation, causeur2011factor}. The method of Friguet {\it et al.} performs a decorrelation of the covariates and compute corrected covariates that are suitable for testing and/or regression. After that, for the selection and the ranking, we propose to combine different selection methods. 

%According to these considerations, we propose in this paper a general methodology in order to select and rank quantitative covariates that are linked to an outcome $Y$. 
Our methodology is divided in two steps: a pretreatment of the covariates (step 1) and a procedure of selection of the pretreated covariates (step 2). The pretreatment consists of (step 1.1) detecting the clusters of covariates by using the clustering of covariates proposed by \cite{chavent2011clustofvar}, and (step 1.2)   applying a "decorrelation" between the covariates inside each cluster using the factor analysis proposed by \cite{friguet2009factor, friguet2011estimation, causeur2011factor}. Their method performs a decorrelation of the covariates and calculates the corrected covariates suitable for statistical testing and/or regression. 

After that pretreatment, we propose a procedure to select and rank the covariates, by combining different selection methods that take into account the nature of the outcome $Y$ (qualitative or quantitative) and the high dimensional context (multiple testing procedures for the tests, penalised regression, ...). We define a score for each covariate, which is defined by the number of selections among all the selection methods involved in this step. This score can be used to classify the covariates like in \cite{su2003rankgene}.

The paper is organized as follows. In Section \ref{sec:method}, we detail the model and explain the principle of the main steps of our methodology: the pretreatment of the covariates and the construction of the covariates scores of selection. Section \ref{sec:simu} is dedicated to simulations studies in order to assess the interest of the proposed pretreatment on one hand and the good working of the whole selection strategy on the other hand. The simulations are performed in the case where the variable of interest is binary. Section \ref{sec:data} is dedicated to two real data analysis: the purpose of the first analysis is to select covariates that are linked with the outcome of a lung cancer treatment, whereas the second analysis selects covariates linked with a breast cancer biomarker. In both analysis, the selected covariates are used to define genetic profiles of patients. Section \ref{sec:conclu} gives some conclusions and perspectives. The Appendix gives two simulation studies in the case where the variable of interest is respectively a binary and a continuous quantitative variable (Sections \ref{sec:classif_Appendix} and \ref{sec:appendix} of Appendix). Technical details on the two real data applications are also given (Sections \ref{SM:transgene} and \ref{SM:bio} of the Appendix).

\section{Methodology}
\label{sec:method}

\subsection{Framework and model}
\label{model}
 We suppose that we have $n$ i.i.d replications of $(Y,{\mathbf X} )$ where $Y$ is the variable of interest, and $ {\mathbf X} =(X_1,X_2,\ldots, X_p) $ is the vector of covariates, taking its values in $\mathbb R^p$. We make the assumption that the covariates are decomposed into $K$ independent clusters: 
 $$\mathbf X= (\underbrace{X^{(1)}_1, \ldots, X^{(1)}_{p_1}}_{\mathbf X^{(1)}}, \ldots, \underbrace{\ldots, X^{(k)}_i,\ldots}_{\mathbf X^{(k)}}, \ldots, \underbrace{\ldots,X^{(K)}_{p_K}}_{\mathbf X^{(K)}}) = (\mathbf X^{(1)}, \ldots, \mathbf X^{(K)}),$$

where $p_1 + \ldots + p_K=p$.

On one hand, we model the dependence in the $K$ clusters of covariates as in the framework of \cite{friguet2009factor}: inside each cluster $\mathbf X^{(k)}$, the common information between the $p_k$ covariates is modeled by regression on $q_k$ latent factors:
\begin{equation}\label{frigu}
X^{(k)}_i = \delta^{(k)}_i(Y)+\mathbf{b}^{(k)}_i \mathbf Z^{(k)}+\varepsilon^{(k)}_i, \qquad \textup{for } i=1, \ldots, p_k, 
\end{equation}
where $\delta^{(k)}_i(Y)$ is a function of $Y$,  $\mathbf Z^{(k)}$ is a random centered $q_k$-vector such that $\mathbb E(\mathbf Z^{(k)} \mathbf Z^{(k)'}) = \mathbf{I}_{q_k}$, $\mathbf{b}^{(k)}_i$ is a $q_k$-vector, and $\mathbf{\varepsilon}^{(k)}=(\varepsilon^{(k)}_1,\ldots,\varepsilon^{(k)}_{p_k})$  is a random centered $p_k$-vector with independent components, and independent of $\mathbf Z^{(k)}$.  The common information contained in $\mathbf X^{(k)}$ is then  concentrated in a small dimension space by $q_k$ latent factors $\mathbf Z^{(k)}$. Under the model (\ref{frigu}), we have:
\begin{eqnarray}\label{Sigma}
\mathbf{\Sigma}^{(k)}:=\mathbb{V}(\mathbf X^{(k)}|Y)  &=&\mathbf{B}^{(k)} (\mathbf{B}^{(k)})'+\mathbf{\Psi}^{(k)}\\
\mathbb{V}(\mathbf X^{(k)}|Y,\mathbf Z^{(k)})  &=&\mathbf{\Psi}^{(k)}\\
\mathbb{C}ov(\varepsilon_i^{(k)},Z_j^{(k)})&=&0,\quad \forall i,j,k
\end{eqnarray} 
where $\mathbf{\Psi}^{(k)}$ is a diagonal $p_k\times p_k$ matrix (the covariance matrix of $\mathbf{\varepsilon}^{(k)}$) and $\mathbf{B}^{(k)}$ is a $p_k\times q_k$ matrix of factor loadings $\mathbf{b}^{(k)}_i$ (cf Equation (\ref{frigu}), the $\mathbf{b}^{(k)}_i$ being the $i$th row of $\mathbf{B}^{(k)}$). In the decomposition given in Equation (\ref{Sigma}), the diagonal element $\Psi^{(k)}_i$ is the specific variance of the response $X^{(k)}_i$ while $\mathbf{B}^{(k)}(\mathbf{B}^{(k)})'$ appears as the shared variance in the common factor structure. \cite{friguet2009factor} define the common variance by 
\begin{equation}\label{eq:cv}
\textup{ComVar}^{(k)}=\frac{\textup{trace}(\mathbf{B}^{(k)}(\mathbf{B}^{(k)})')}{\textup{trace}(\mathbf{\Sigma}^{(k)})}.
\end{equation} 

\begin{figure}[!p]
	\centering
	\includegraphics[scale=0.5]{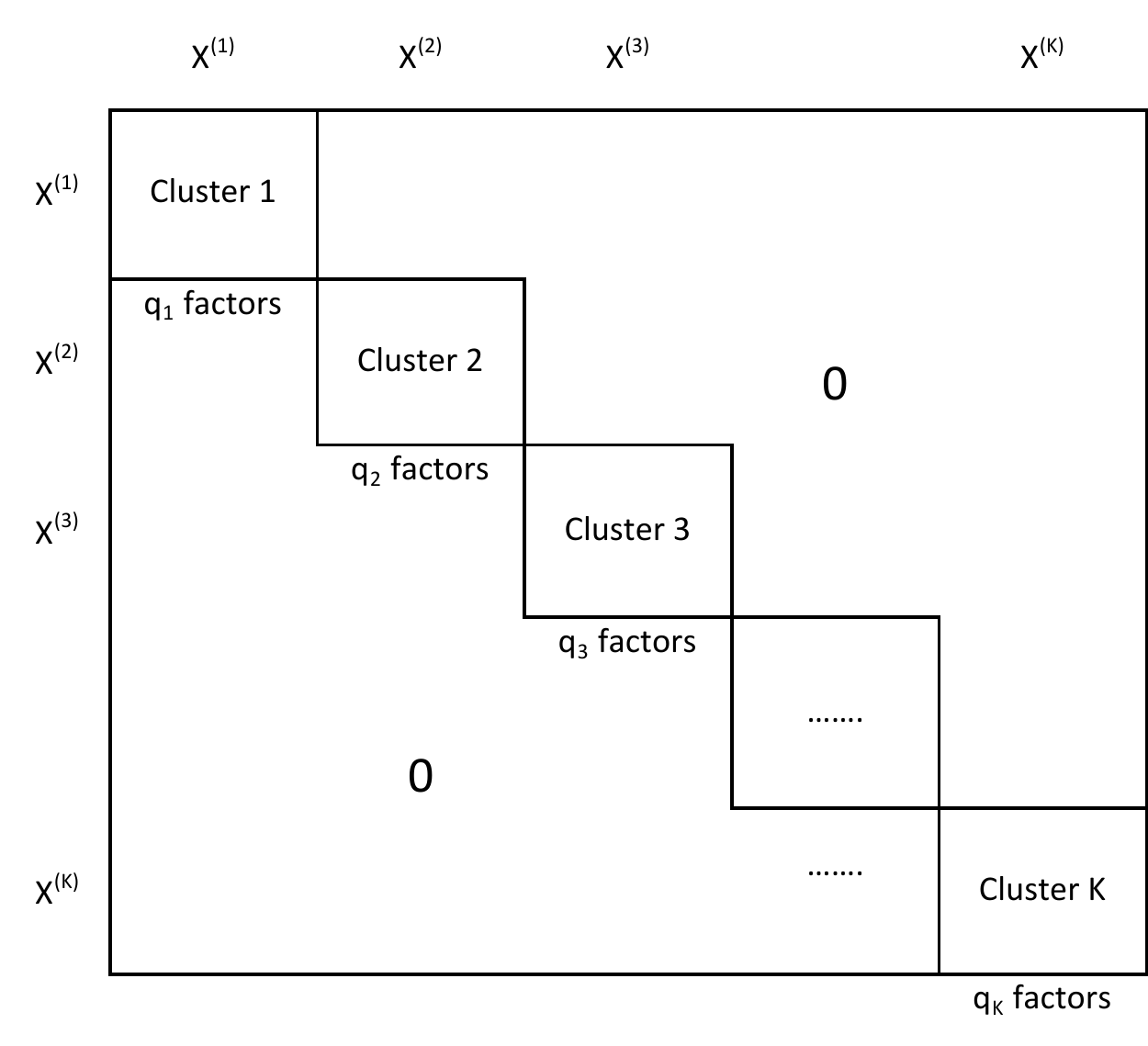}
	\caption{Covariance structure of covariates given $Y$}
	\label{fig:correl}
\end{figure}

On the other hand,  we suppose that the specific informations  at each cluster (that is vectors $(\mathbf Z^{(k)},\mathbf{\varepsilon}^{(k)})_{1\leq k \leq K}$) are independent, then, given $Y$, the covariance matrix of the whole vector of covariates has the form given by the Figure \ref{fig:correl}. 

\subsection{Main prodecure}

The procedure is decomposed in two steps: a pretreatment of the covariates (step 1) and a selection method of the covariates (step 2). 

%\subsubsection{Step 1: Pre-filtering of covariates}
% The dimension is extremely high, meaning very large $p$. It is necessary to pre-filter in order to reduce the dimension as a first step. The dimension reduction is necessary for two reasons : firstly, it is usual to assume that the number of significant covariates is sparse, and secondly, computational cost is always the biggest hindrance in high-dimensional data analysis. From a computer science point of view, when the dimension is extremely high, it is actually a NP-hardness, which cannot be solved in polynomial time. Many existing methods and R packages, e.g. clustering, random forest, etc., are only capable of dealing with a relative high dimension. In the studies we have in mind, the data have extreme high dimension ($p$ is in the order of tens of thousand) with a limited sample size ($n$ between 25 and 100). Algorithms, like lasso, partial least squares etc., are recommended, since they assume a sparse scenario which notably reduces the computational complexity. However, those algorithms have certain limits in high dimension, which will be discussed in section \ref{sec:simu}. 

\subsubsection{Step 1: pretreatment of data (clustering of covariates and decorrelation inside clusters).}

The aim of this pretreatment is to perform a decorrelation of the covariates, to obtain corrected covariates that are suitable for testing and/or regression.  Indeed, the correlation between covariates has an impact on all the classical selection procedures: the conventional methods, namely the multiple testing procedures (the p-value adjustment methods proposed by \cite{bonferroni1936teoria}, \cite{benjamini1995controlling,benjamini2001control}, or the q-value proposed by \cite{storey2002direct,storey2003statistical,storey2003positive}, or the local FDR presented in \cite{aubert2004determination}, \cite{bar2005comparaisons}) are all built on the assumption that tests are independent. As a results, they are no longer promising if the independence is not verified. A very detailed discussion can be found in the Friguet's thesis \cite{Friguet2012}.

In estimating together the latent factors $\mathbf Z^{(k)}$ and the coefficients of regressions ($\mathbf{B}^{(k)},\mathbf{\Psi}^{(k)}$) by an E.M. algorithm in model  (\ref{frigu}), the FAMT procedure of \cite{friguet2009factor}  can correct the covariates such that they are almost independent and as a result, suitable for multiple testing procedures or selection by regression or random forests. More precisely, the corrected data, noted $X^{(k)\star} _i=X^{(k)} _i-\mathbf{b}^{(k)}_i{\mathbf Z}^{(k)} = \delta^{(k)}_i(Y)+\varepsilon^{(k)}_i$, $i=1,\ldots,p_k$, lead to a standard multiple regression problem where the errors $\varepsilon^{(k)}_i$ are independent. Note that this correction  of the data ${\mathbf X}$ is done conditionally on the variable of interest $Y$ \cite{friguet2009factor}.

Of course, the whole vector $\mathbf X$ satisfies assumption of Equation (\ref{frigu}), and \cite{friguet2009factor}   apply this decorrelation procedure on the whole set of covariates $\mathbf X$. But instead of applying Friguet's procedure on the whole set of covariates $\mathbf X$, we propose to first detect the different clusters $(\mathbf X^{(k)})_{1 \leq k \leq K}$ and then to apply the decorrelation method on each cluster. Indeed,  \cite{yaojie} has shown with some simulation studies that the decorrelation was degraded by the dimension of the vector of covariates, whereas it was better after the detection of the independent clusters. By this way, the covariates selection procedure can be highly improved  by clustering of covariates (step 1.1) before applying factor analysis to correct the correlation within each cluster (step 1.2), as it is shown in Section \ref{sec:simu}.

\paragraph{Step 1.1: clustering of covariates}
We apply a clustering of covariates in the purpose to find clusters of correlated variables as we assumed in Section \ref{model}. We propose to use the  algorithm of \cite{chavent2011clustofvar} to cluster covariates into homogeneous clusters and thus to reveal structures. This algorithm maximizes an homogeneity criterion, where the homogeneity of a cluster is defined by the sum of squared Pearson correlations between the covariates present in the cluster and the first principal component of this cluster. This algorithm is expected to roughly find the highly correlated clusters of covariates as we assumed in the Section \ref{model}. The procedure proposes also a method (based on bootstrap resampling) to find the number $K$ of clusters if it is unknown.

\paragraph{Step 1.2: Factor analysis to correct dependency structure in each cluster}

As already explained in the beginning of this section, clustering is followed by decorrelation inside each cluster using the Friguet's procedure.

At the end of this pretreatment procedure, we obtain corrected data, noted ${\mathbf X}^*_Y$ in the sequel. Note that ${\mathbf X}^*_Y$ depends on $Y$.  To simplify the notations, ${\mathbf X}^*_Y$ will be noted ${\mathbf X}^*$.

\subsubsection{Step 2: Aggregation of statistical methods applied on the resulting dataset.}\label{Methode}

The statistical methods proposed in this part are not fixed and can be adapted by the practitioner according to its preferred selection methods and the characteristics of the data (nature of variable of interest $Y$, samples' sizes and so on).

The idea is the following: we choose several methods to select the pretreated covariates $\mathbf X^{\star}$. We perform $L$ methods, then each covariate $X^{\star}_i$ obtains a score $S_i \in \{0,1, \ldots , L\}$ that is the number of selections among the $L$ methods. By this way,  the covariates can be ranked according to their link with the outcome $Y$.

For instance, in the examples developed in our simulation studies and in real data, $Y$ is binary, the size of the samples are low and we choose eight different methods of selection: five different multiple testing procedures applied to the Wilcoxon test (Bonferroni, Benjamin-Hochberg, q-values, local FDR, FAMT), logistic regression penalised by Lasso, and two selections by random forests (threshold step and interpret step, see \cite{genuer2010variable}).
 The outcomes of this procedure are the scores  $S_i, i=1,\ldots,p$ which are integers included in $[0,8]$. For example, if $S_i=8$, then the corresponding variable has been selected by all the eight methods, whereas if $S_i=0$, the corresponding variable has been selected by none of them. The scores can be used to rank the covariates according to the strength of their link with the variable of interest.  

\subsection{R Package \texttt{armada}}
In the sequel, we call our procedure ARMADA for AggRegated Methods for covAriates selection under Dependence. Our procedure has been implemented in an R package, called \texttt{armada}, available on the CRAN \cite{package_armada}. The package proposes also a graphical representation of the selected covariates, through heatmaps, as presented in Figures \ref{fig:heatmap_score_reg_classif2} and \ref{Fig3HD}.

\section{Simulations}
\label{sec:simu}

We first explain the simulation design in Section \ref{sec:designs}. We then describe the effect of the pretreatment in Section \ref{sec:pretreatment} and finally, we study the selection procedure in Section \ref{sec:results}.

\subsection{Simulation design}\label{sec:designs}
We propose a simulation study with $p=1600$ covariates and sample size $n=60$. We first describe how to create dependence in the covariates $\mathbf X$, then we present a simulation design in a classification study. Two other designs in classification and regression cases are given in Sections \ref{sec:classif_Appendix} and \ref{sec:appendix} of the Appendix.

The covariates $\mathbf X=(\mathbf X^{(k)})_{k=1,\ldots,4}$ are clustered into four independent clusters, each of them containing $p_k=400$ covariates. For this, before to model the dependence with the outcome $Y$, we generate for each cluster $k$,  a preliminary vector $\mathbf{ \tilde{X}}^{(k)}$ that is a gaussian $400$-vector, with mean $0$ and non-diagonal variance-covariance matrix  $\mathbf{\Sigma}^{(k)}$. The correlation between the covariates of $\mathbf{\tilde{X}}^{(k)}$ inside the cluster $k$ is designed by a factor analysis model described in Equation (\ref{Sigma}). More precisions on the simulation procedure of data with covariance design defined by Equation (\ref{Sigma}) can be found in  \cite{Friguet2012}. We simulate data with common variances  $\textup{ComVar}^{(k)}$ equal to $0.8$ in each cluster (recall that the common variance is defined in Equation (\ref{eq:cv})).  Moreover, the numbers of latent factors in each cluster are $(q^{(1)},\ldots ,q^{(4)})=(4,6,8,10)$.

Now, we create the dependence between $\mathbf X$ and $Y$ in perturbing some component of $\mathbf{ \tilde{X}}$. We consider an equiprobable two-class problem, $Y\in \{0,1\}$ (i.e. $Y=1$ for $\frac{n}{2}$ subjects, and $Y=0$ for $\frac{n}{2}$ subjects). $Y$ is linked with 160  influential covariates in $\mathbf X$, whose links with the response variable $Y$ have different intensities. The other covariates are noise. 
More precisely,  in each cluster $k=1,\ldots, 4$, and for $i=1,\ldots, 400$, $${X}^{(k)}_i = \tilde{X}^{(k)}_i+\delta_i\mathbf{1}_{Y=0}$$
with $\delta=(\delta_i)_{i=1,\ldots,  400} = (\underbrace{1.5, \ldots, 1.5}_{i=1,\ldots,10}, \underbrace{1, \ldots, 1}_{i=11,\ldots,20},\underbrace{0.75, \ldots, 0.75}_{i=21,\ldots,30},\underbrace{0.5, \ldots, 0.5}_{i=31,\ldots,40}, \underbrace{0, \ldots, 0}_{i=41,\ldots,400})$. In other words, $Y$ is linked with the $m_1=40$ first covariates of each cluster, and  the $m_0 = 360$ remaining covariates of each cluster are independent of $Y$. Then, the $10$ first covariates of each cluster are the most strongly linked with  the response variable $Y$  and the  strength of the link is decreasing in the successive groups of $10$ influential covariates.

We can remark that this design respects the covariance matrix given in Figure \ref{fig:correl} and the  model given by Equation (\ref{frigu}).

%%%%%%%%%%%%%%%%%%%%%%%%%%%%%
\subsection{Interest of our data pretreatment}\label{sec:pretreatment}

\noindent
In order to emphasize the interest of our data pretreatment, we compare the results of a Wilcoxon test after three different  data pretreatments:

\begin{enumerate}
\item[Procedure 1:] nothing is done on the dataset $\mathbf X$.
\item[Procedure 2:]  the  covariates $\mathbf X$ are decorrelated, taking $Y$ into account, with the factor analysis procedure of  \cite{friguet2009factor,causeur2011factor}, implemented in the R package  \texttt{FAMT}. This gives a new dataset $\mathbf X^\dagger_Y$.
\item[Procedure 3:] the 4 clusters are estimated with the procedure of \cite{chavent2011clustofvar}, implemented in the R package  \texttt{ClustOfVar};  then the covariates are decorrelated in each cluster,  taking $Y$ into account, with the factor analysis procedure of \cite{friguet2009factor,causeur2011factor}, implemented in the R package  \texttt{FAMT}. This gives a new dataset ${\mathbf X^*}_Y$ obtained by the  concatenation of the decorrelated clusters.
\end{enumerate}

\textit{\textbf{Remark:} our data pretreatment is the Procedure 3. We have supposed that the number of clusters is known. If that is not the case, the user can choose its own number of clusters by  using the graphical tools of the \texttt{ClustOfVar} procedure (plots of the dendrogram).}

\begin{figure}[htbp]
	\centering
	\includegraphics[scale=0.3]{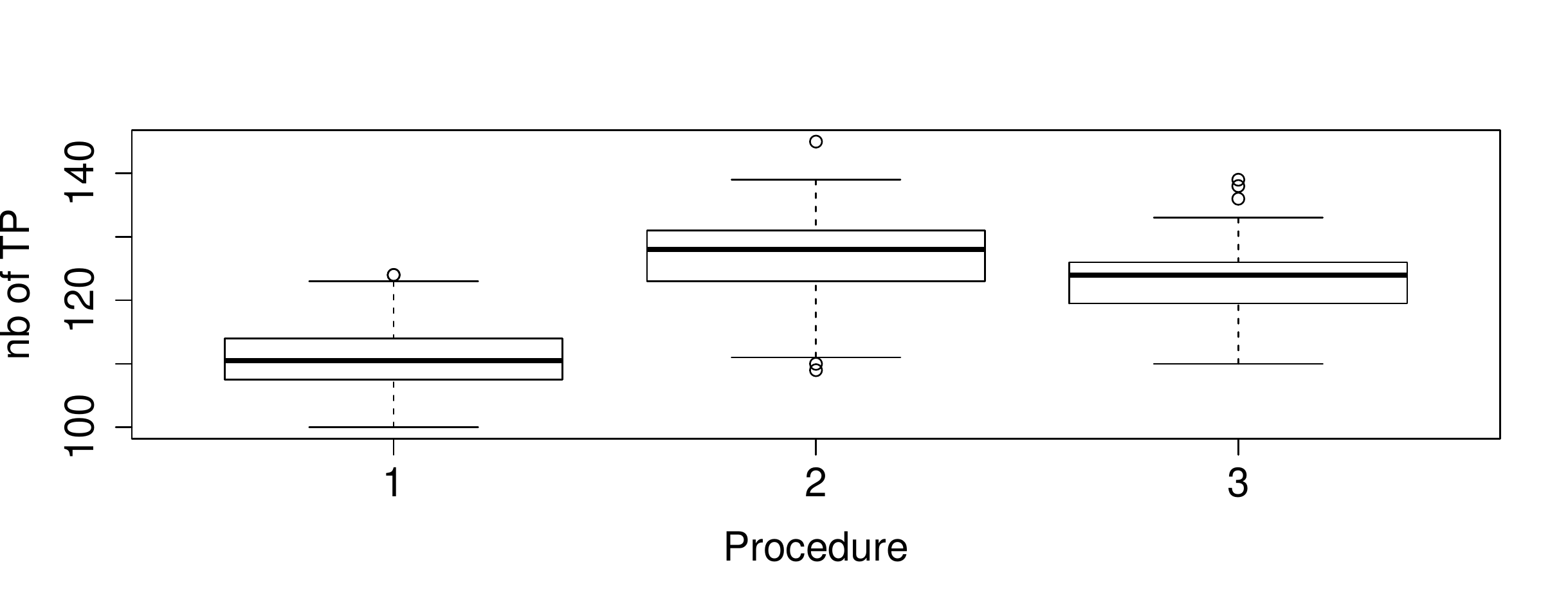}
		\includegraphics[scale=0.3]{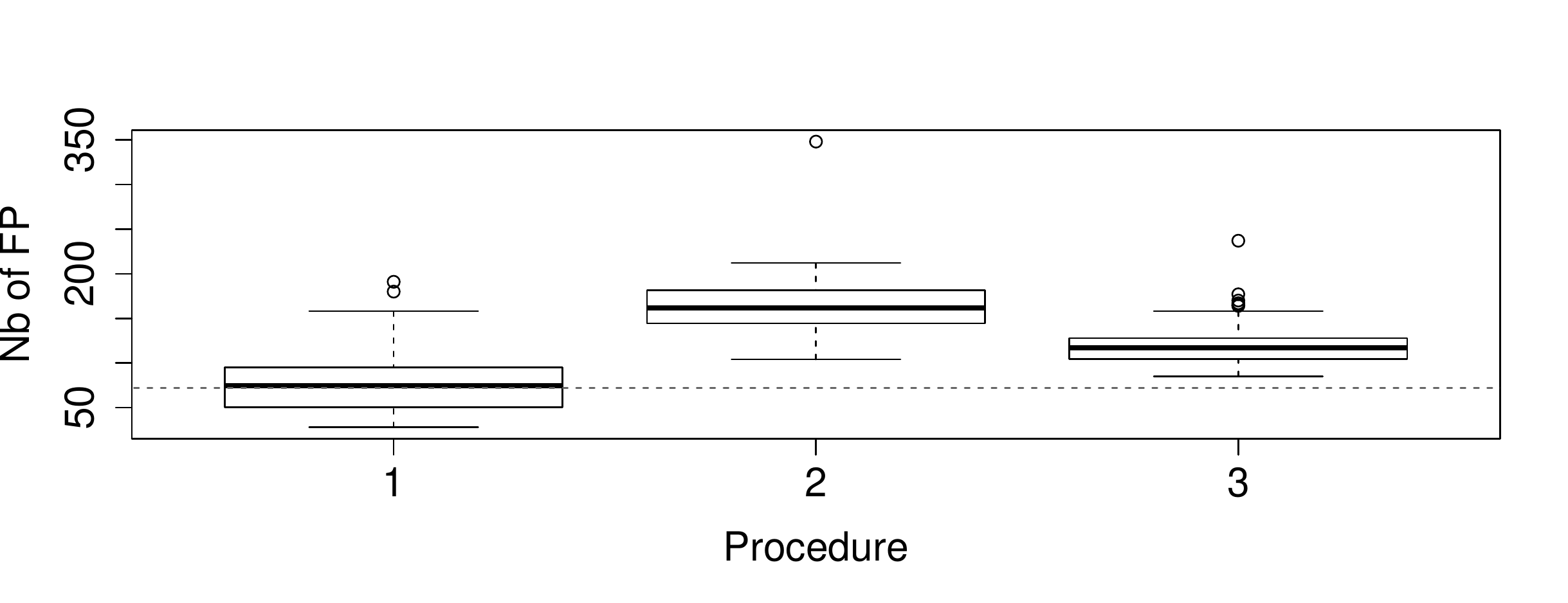}
		\caption{Number of true positive tests (left), false positive tests (right)  according to the different pretreatment procedures (1: Nothing, 2: FAMT, 3: clustering followed by FAMT in each cluster). Dotted lines: expected number of FP. Boxplots are calculated on $N=100$ runs.}
	\label{fig:nb_Positifs_Hafid}
\end{figure}

Our objective is to find out the differently expressed covariates in the two groups (groups $Y=0$ and $Y=1$) with sample sizes $\frac{n}{2}=30$. For this, we perform Wilcoxon tests  on each of the $p$ pretreated covariates of the dataset (that is $\mathbf X$ for Procedure 1, $\mathbf X^\dagger_Y$ for Procedure 2, ${\mathbf X^*}_Y$ for Procedure 3), given a three sets of $p$ p-values. For each of these procedures, the selected covariates are those with p-values lower than  0.05. We compare these procedures on $N=100$ runs of $(\mathbf X,Y)$. For the comparison, we count the number of influential covariates that are correctly detected (this number is noted TP, for True Positive), this indicator gives an idea of the sensibility of the test after the procedure. To assess the specificity, we count the number of non-influential detected covariates (this number is noted FP, for False Positive). Note that the perfect method would detect all the influential covariates (that is 160 here) and no False Positive. However, according to the detection threshold chosen for the p-value, the expected number of FP is $72=5\%\times(1600-160)$. The results are shown in Figure \ref{fig:nb_Positifs_Hafid}.

On the Figure \ref{fig:nb_Positifs_Hafid}, we can see that Procedure 1 is in fact the one that has the lowest rate of FP but its power is also the poorest. Our Procedure reduces the mean and the variability of the distributions of the false positive rates, in comparison to the Procedure 2 (i.e. the FAMT procedure). The power of our Procedure is comparable with Procedure 2. This results show the interest of our proposed pretreatment before performing selection.

%%%%%%%%%%%%%%%%%%%%%%%%%%%%%
\subsection{Results of the whole method (pretreatment and selection)}\label{sec:results}

\begin{figure}[htbp]
	\centering
		\includegraphics[scale=0.3]{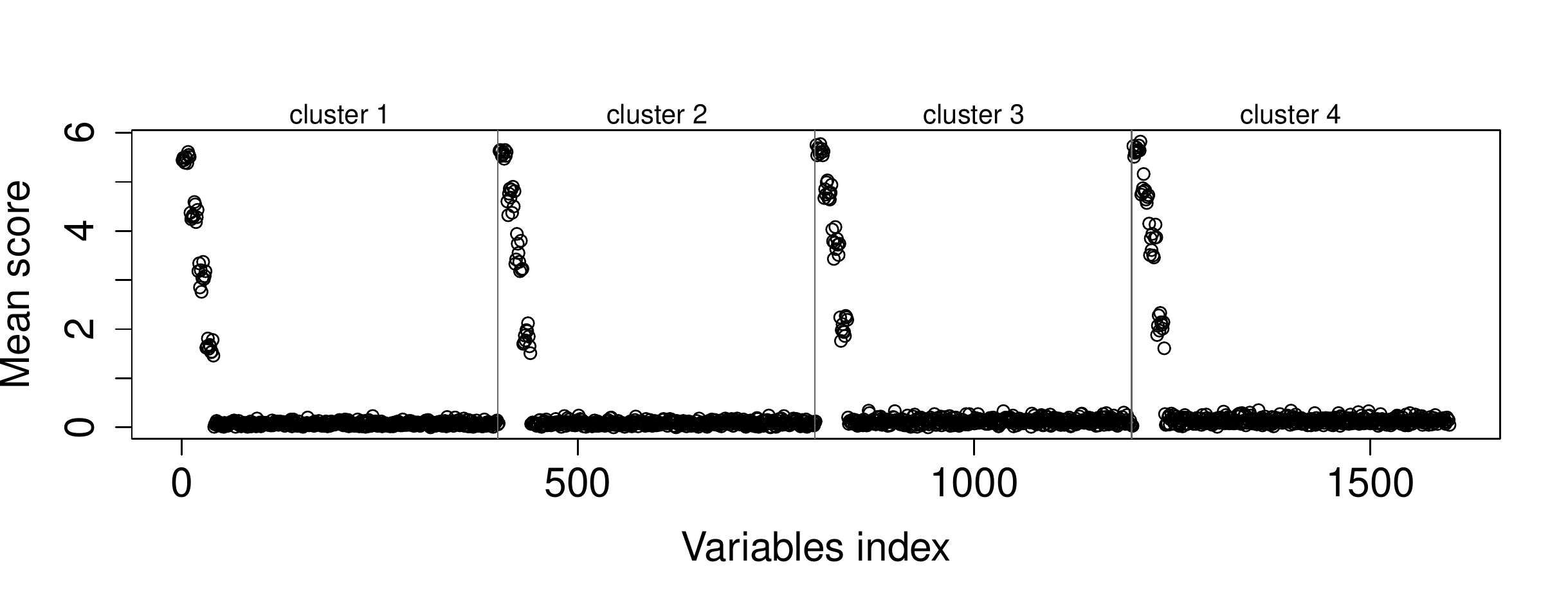}
				\includegraphics[scale=0.3]{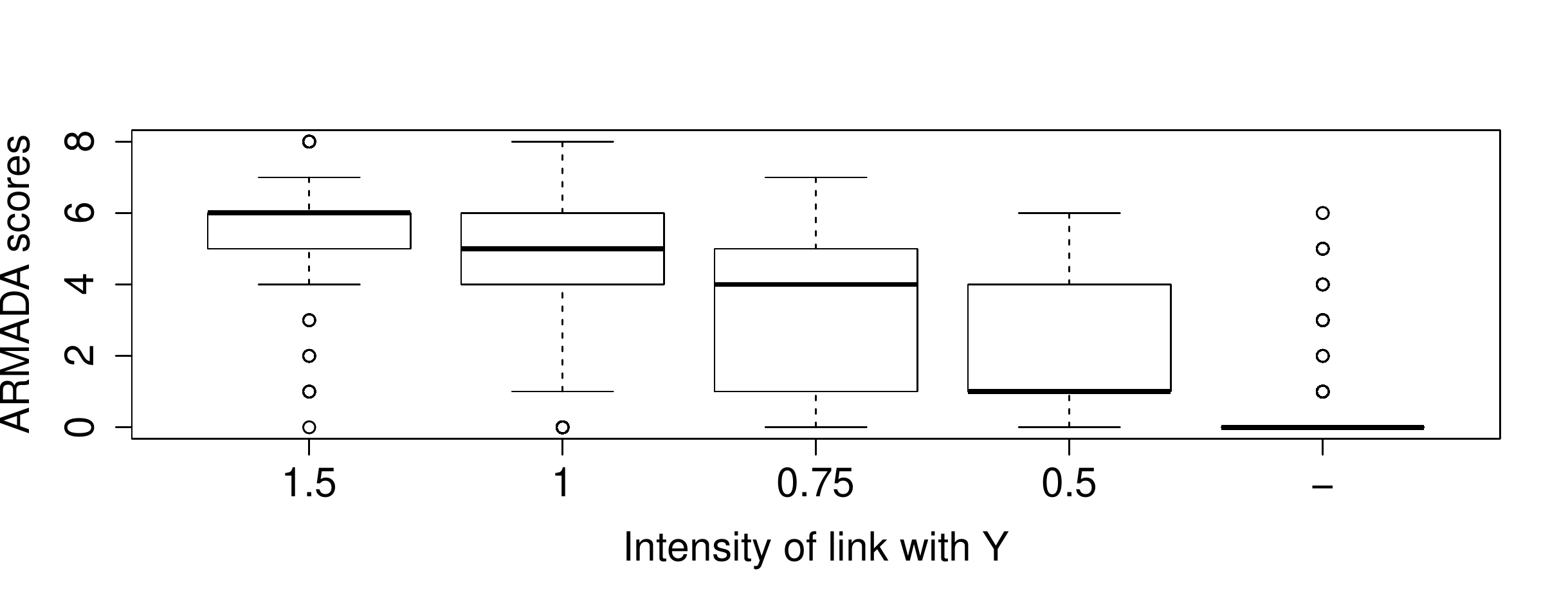}
			\caption{Left: mean of the ARMADA scores obtained by all the covariates. Right: boxplot of the scores of the covariates, ranked by levels of link with $Y$. Means and boxplots are calculated on $N=100$ runs.  }	\label{fig:mean_score_Hafid}
\end{figure}

In order to describe the performances of our method, we show in Figure \ref{fig:mean_score_Hafid} the mean ARMADA scores obtained on the $N=100$ runs of $(\mathbf X, Y)$ for each design. The  scores are given for  all the covariates individually, and also by group of influential and noise covariates (the groups of influential covariates are noted by "1.5", "1", "0.75", "0.5" (see Section \ref{sec:designs});  the group of noise covariates is noted by "-").

We can see on the Figure \ref{fig:mean_score_Hafid} that the scores give a clear ranking of the covariates, according to the strength of their link with the response variable $Y$. The highest scores are obtained by the covariates which are the most strongly linked with the response variable $Y$. The ARMADA method is  performant: the mean score clearly distinguishes the five groups of covariates according to their link with $Y$. The distribution of the individual scores inside each group is given by the boxplots. The ARMADA scores clearly separate the influential covariates from the others; and inside the influential covariates the two first groups are clearly separated of the last one. 
Note that around 95\% of the noise covariates obtained an ARMADA score that was exactly 0.

\subsection{Comparison with other selection methods}\label{sec:comparison}
We propose the following selection criterion in our procedure: the selected covariates are those with scores greater or equal to 1. 

We compare this selection procedure with two other selection methods: 
\begin{itemize}
\item the Wilcoxon test: the selected covariates are those with raw-pvalues (i.e. p-values without any correction) lower than 0.05,
\item the FAMT procedure \cite{causeur2011factor}: the selected covariates are those with adjusted p-values lower than 0.05. 
\end{itemize}

To compare the three selection methods, the Table \ref{tab:recap2_Hafid} gives the rates of selection for each group of influential covariates, and for the group of noise covariates.  The rates of selection have been computed on $N=100$ runs of $(\mathbf X,Y)$. We can see that our method respects the expected rate of false positives that is not the case for the FAMT method which exhibits a greater rate of 10 \%.  Moreover, our method gives the best results.  The rate of selection of the influential covariates is very good compared with the other methods even if the strength of the link is poor.

\begin{table}[htbp]
\tbl{Results of the $N=100$ runs: rates of selection of the different groups of influential and noise covariates by the ARMADA method, the Wilcoxon test and the FAMT procedure. The corresponding standard deviations are given in brackets. }
{\begin{tabular}{cccc}\toprule
 & ARMADA & Wilcoxon & FAMT\\ \midrule
 1.5&                      0.99 (0.04) &                     0.99 (0.07)&                       0.99 (0.02)\\
 1&                        0.97 (0.15) & 0.85 (0.35) & 0.95 (0.20)\\
 0.75&                     0.91 (0.27) &                     0.62 (0.48) &                       0.82 (0.38)\\
0.5&                      0.79 (0.40) &                     0.33 (0.47) &                       0.52 (0.49)\\
 -&                        0.05  (0.23)&                     0.05 (0.22)&                       0.10 (0.30)\\ \bottomrule
\end{tabular}}
\label{tab:recap2_Hafid}
\end{table}

Finally, we can conclude with the ROC curves given in Figure \ref{fig:ROC}  that our method outperforms the two others selection methods (the ordinates of the points of the ARMADA ROC curve are all higher than the ordinates of the points of the two other ROC curves). The ROC curves have been obtained by the mean of the $N=100$ ROC curves obtained in the $N=100$ runs of $(\mathbf X,Y)$.

 \begin{figure}[htbp]
	\centering
		\includegraphics[scale=0.4]{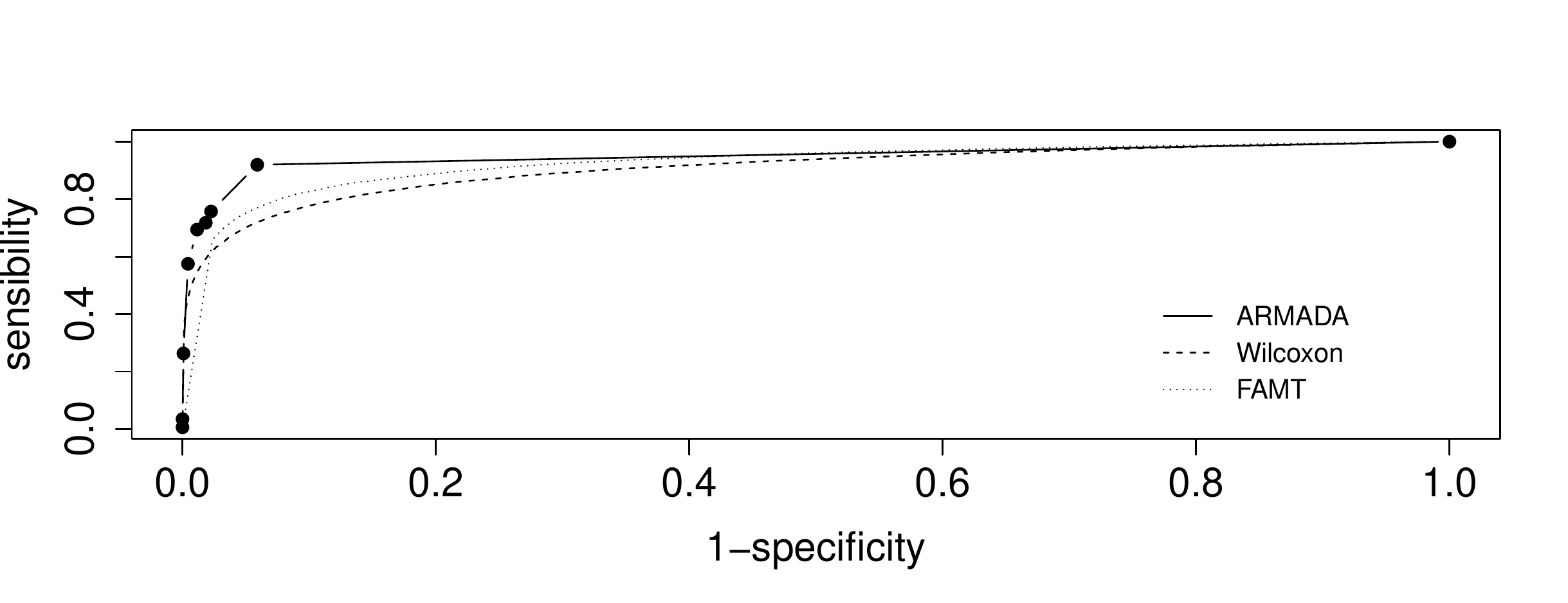}
			\caption{ROC curves for the three selection methods. }	\label{fig:ROC}
\end{figure}

\section{Application to real data}
\label{sec:data}

In this section we apply  our method on two real datasets, both in oncology.  The first one concerns the selection of transcriptomic covariates linked with  the outcome of a chemotherapy for lung cancer. The second one concerns the selection of covariates linked with a quantitative biomarker ER$\alpha36$ in breast cancer.

\subsection{Outcome of chemotherapy for advanced non-small-cell lung cancer}
\label{sec:data1}
 We apply our method on transcriptomic data of $n=37$ patients with advanced non-small-cell lung cancer, who have received chemotherapy.  Even if we are aware of the fact that chemotherapy is not a target therapy, the problematic is really to select suitable transcriptomic covariates in the purpose to detect profiles associated with the effect of a treatment.  For each patient, we have 51336 transcriptomic covariates, and its survival status: 24 patients whose death occurred before 12 months  and 13 patients whose  death occurred after 12 months. This criteria of death before one year is very common in clinical trials. We applied a first filtering of the covariates, where we decided to ignore the covariates for which the Wilcoxon test does not detect a difference between the 24 patients whose survival time is lower than 12 months and the 13 other ones (we eliminate covariates with Wilcoxon-pvalue greater than 0.05). After this filtering we obtained a dataset with $n=37$ patients and $p=6810$ covariates. In the pretreatment step, we found  that the $p=6810$ covariates are decomposed in 3 independent group of covariates.

In a first time, the biological question was to find the genes which can explain a survival time greater or lower than 12 months. We then consider a binary response variable $Y$:  $Y=1$ for the 24 patients whose survival time is lower than 12 months and $Y=0$ for the 13 patients whose survival time is greater than 12 months. The response variable being binary, we first applied the selection method presented in Section \ref{Methode}.  %The results are shown in Table \ref{tab:res_transgene}. We see that 10 covariates are particularly important, with a score equal to 7, whereas 2827 covariates have a score equal to 0, and 3983 covariates have a score greater or equal to 1.

%It is clear that, the biologist will not focus on the 3983 covariates with a positive score. But the method clearly gives a hierarchy between the genes and it is sure that the function of the 10 genes with a score at 7 has to be studied to understand its link with the "success" of the treatment. 

Moreover, as the survival time was known for all the 37 patients without any censoring, we also apply our method on the same dataset (6810 covariates) but here,  $Y$ is the survival time. We then have a regression problem. We have used eight selection methods in Step 2 of our method: five different multiple testing procedures applied to the Pearson correlation test (Bonferroni, Benjamin-Hochberg, q-values, local FDR, FAMT),  regression penalised by Lasso, and two selections by random forests (threshold step and interpret step, see \cite{genuer2010variable}).  %The results are given in Table \ref{tab:res_transgene_reg}.

\begin{figure}[htbp]
	\centering
		\includegraphics[scale=0.5]{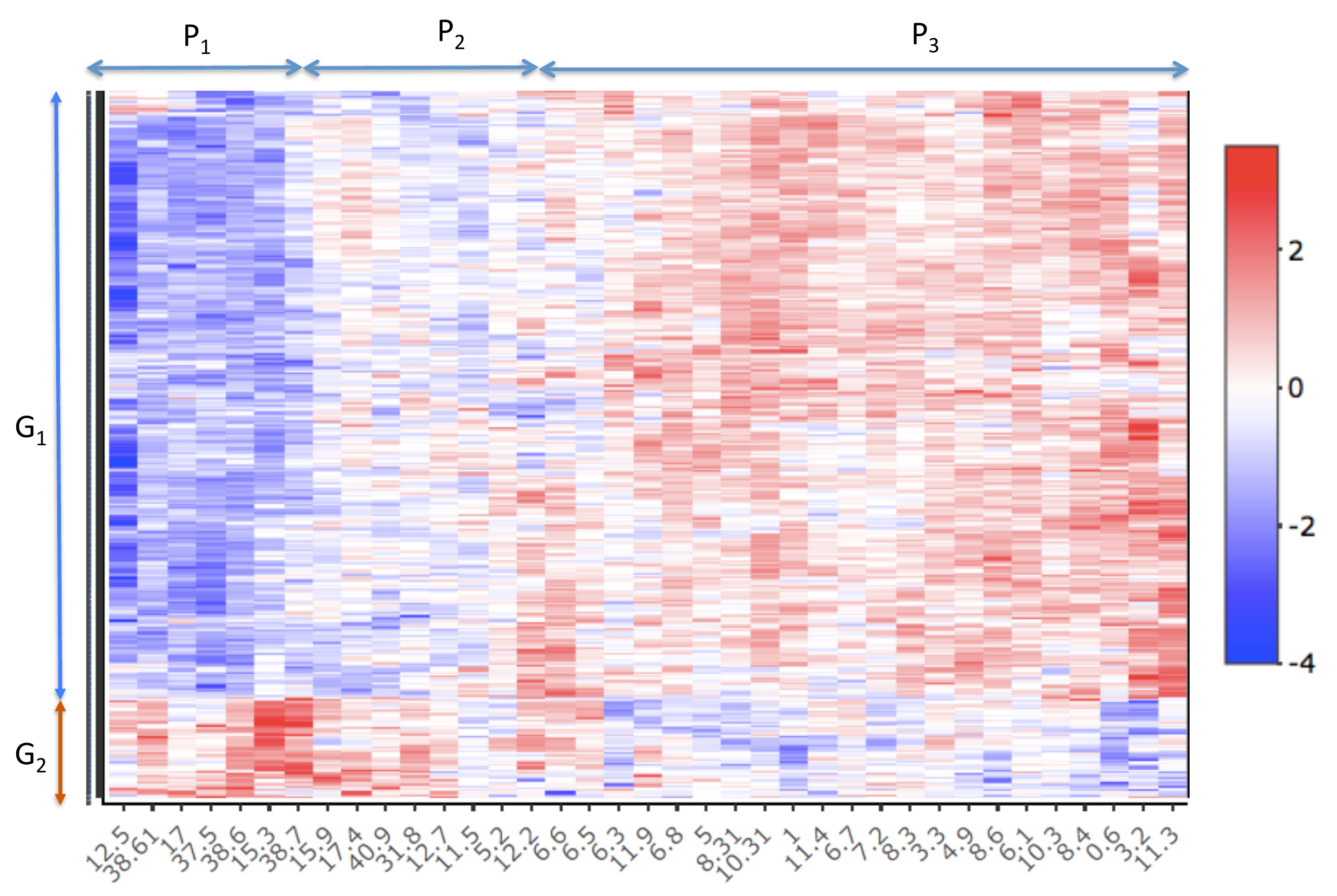}
		\caption{Heatmap of the 342 covariates which have ARMADA scores greater or equal to 5 in \textbf{both} classification \textbf{and} regression studies. Each column corresponds to one patient. The $x$-axis represents the patients (marked with their survival time) and the $y$-axis the covariates. Color gives the normalised expression of the covariates.}
	\label{fig:heatmap_score_reg_classif2}
\end{figure}

The joint results of  the classification and regression studies are given in Table \ref{tab:res_transgene_classif_reg}. In the classification study,  we can see that 10 covariates are particularly important, with a score equal to 7, whereas 2827 covariates have a score equal to 0, and 3983 covariates have a score greater or equal to 1. It is clear that, the biologist will not focus on the 3983 covariates with a positive score. But the method clearly gives a hierarchy between the genes and it is sure that the function of the 10 genes with a score at 7 has to be studied to understand its link with the "success" of the treatment. The table \ref{tab:res_transgene_classif_reg} is a little disappointing, because regression and classification do not select the same covariates. Whatever, among the covariates with a C-score (score in the classification case) equal to 7, there is only one with a R-score (score in the regression case) lower than 4 (equal to 0!). But these two analyses are not looking for the same kind of link with the covariates. 
Moreover,  these two approaches give two tools to detect influential covariates. We can combine these two approaches and consider the covariates that are selected by at least one approach, or consider the covariates that are selected by both of them. In the Figure \ref{fig:heatmap_score_reg_classif2}, we show the heatmap of the selected covariates which have a classification score \textbf{and} a regression score greater than five. For the visualisation of the results, we then build an heatmap  obtained thanks to the R package \texttt{heatmaply} after co-clustering of the survival times (on the $x$-axis) and of the covariates (on the $y$-axis) with the function \texttt{hclust} (Figure \ref{fig:heatmap_score_reg_classif2}).

\begin{table}[htbp]
\tbl{Repartition of the covariates scores in the transcriptomic dataset. The R-scores are given in the 9 rows, the C-scores are given in the 8 columns. For instance, 41 covariates have a R-score equal to 1, and a C-score equal to 0.}
{\begin{tabular}{c|ccccccccc}\toprule
 & \multicolumn{8}{c}{Classification score }\\
 Regression score&   0&    1 &   2&    3&    4&    5&    6&    7\\ \midrule
0& 2227&  328&  273&  337&  531&  257&   34&    1\\
  1&   41&    7&    3&    9&   17&   10&   2&    0\\
  2&  131&   35&   39&   52&  119&   71&    9&    0\\
  3  &119&   48&   44&   50&  117&  114&   17&    0\\
  4&  174&   65&   56&   86&  256&  241&  102&    4\\
  5  &119&   64&   40&   57&  116&  176&  116&    4\\
  6   &15&    4&    4&    5&   12&  19&   26&    1\\
  7&    1&    2&    1&    0&    1&    0&    0&    0\\
  8&    0&    0&    0&    0&    1&    0&    0&    0\\ \bottomrule
  \end{tabular}}
  \label{tab:res_transgene_classif_reg}
\end{table}

The visualisation of the co-clustering of the selected genes and the survival leads to the distinction of three different groups of patients (noted P$_1$, P$_2$, P$_3$ in Figure \ref{fig:heatmap_score_reg_classif2}) of respective sizes 7, 8, 22 from the left to the right of the $x$-axis. The co-clustering  identifies also  two clusters of genes (noted G$_1$ and G$_2$ for simplicity).  
All the people except 2 of the two first group P$_1$ and P$_2$ have a life status $Y=1$ (among the two exceptions, one is at the threshold with a survival of 11.5 months), all of the people of the third group P$_3$ have a life status $Y=0$. The selected covariates clearly discriminates groups P$_1$ and  P$_3$. Indeed, the patients of the group P$_1$ have a low expression of the covariates in G$_1$ and a high expression of the covariates in G$_2$  and the inverse for group P$_3$. Patients of group P$_2$ have intermediate expressions according the two others groups.

As the number of patients $n=37$ is small compared to the number of covariates even after filtering ($p=6810$), we have checked our results with a bootstrap study. The results (reported in the Section \ref{SM:transgene} of the Appendix) show that our method is robust:  the distributions of the bootstrapped scores faithfully reproduce the original scores.

\subsection{Biological network involving ER$\alpha36$ in breast cancer}
\label{sec:data2}

 ER$\alpha36$ is a variant of the oestrogen receptor $\alpha$ encoded by the ESR1 locus and  expressed only in humans \cite{wang2005identification}.  ER$\alpha36$ expression and activity have been mainly studied in vitro and in vivo in the context of breast cancer. However, due to the lack of comprehensive transcriptomic data that include ER$\alpha36$, only sparse information is available on factors that could act up- and downstream ER$\alpha36$ in biological networks.
Our challenge from a statistical point of view was to explain the ER$\alpha36$ expression variation obtained in a small number of breast tumors from a large number of potential explanatory variables that correspond to the 54676 transcriptomic probes. For this, we analysed the biological network involving ER$\alpha36$ through the use of 4 sets of Affymetrix transcriptomic data obtained from breast tumors of different molecular subtypes: the triple negative (noted TN), ERa66+, PR+ and PR- datasets  (details are provided in Section \ref{SM:bio} of the Appendix).

\begin{figure}[htbp]
	\centering
		\includegraphics[scale=0.5]{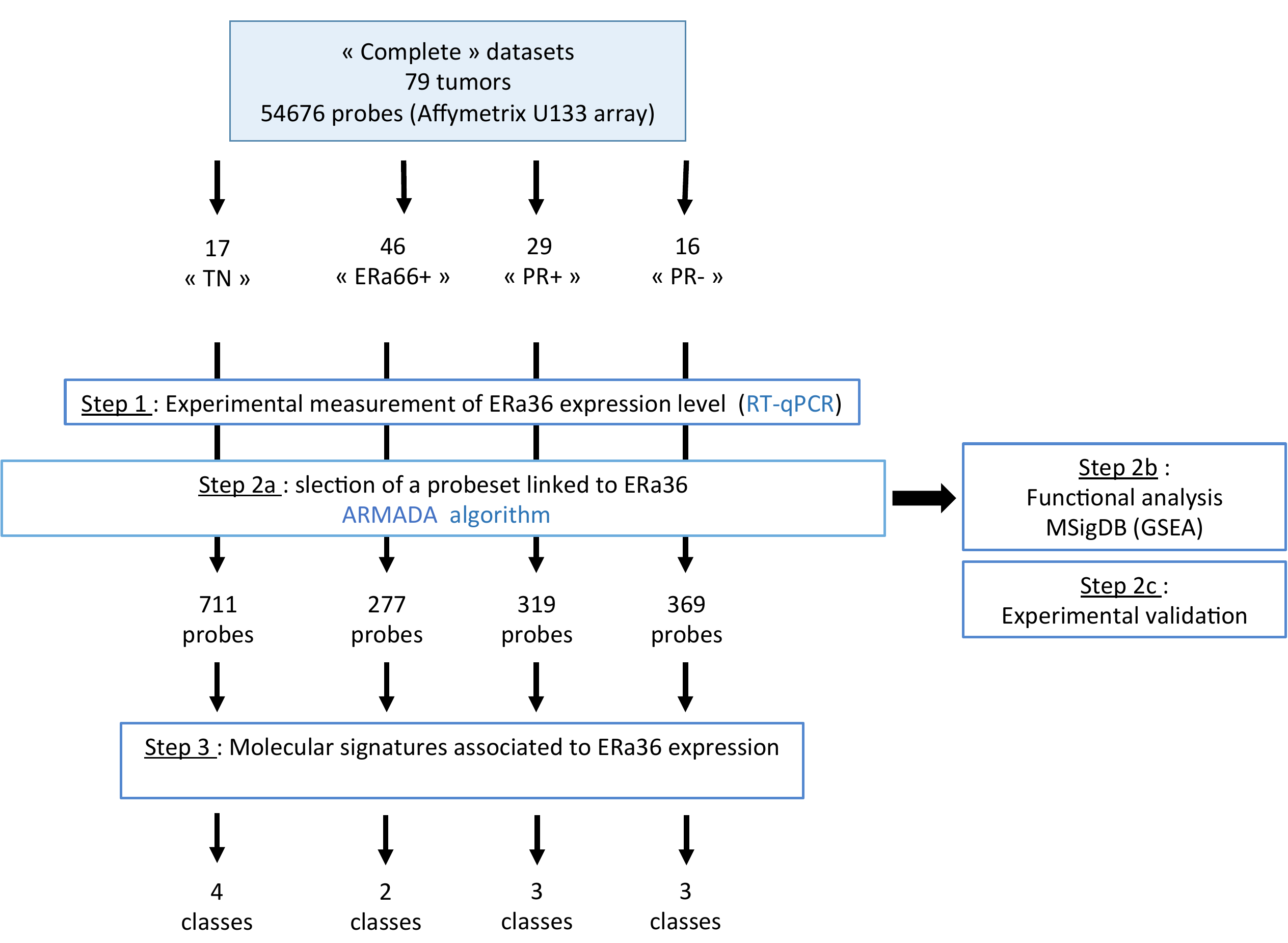}
		\caption{Workflow describing the 3 steps of four breast tumor transcriptomic datasets.}
	\label{Fig1HD}
\end{figure}

The analysis was performed in three steps  (Figure \ref{Fig1HD}).

\noindent
{\bf Step 1, clinical data completion.} These 4 transcriptomic datasets were completed by the measurement of ER$\alpha36$ expression level in each tumor. Biological details are given in Section \ref{SM:bio} of the Appendix.

\noindent
{\bf Step 2a, statistical analysis.} To explain the ER$\alpha36$ expression variation obtained in a small number of tumors from a large number of potential explanatory variables, we used the R package \texttt{armada}, in its regression version (as in Section \ref{sec:data1} where $Y$ is  quantitative), that allowed to select, among the 54676 initial probes, a few hundreds of genes whose expression is supposed to be correlated to that of ER$\alpha36$ (ARMADA score $\geq$ 1). We obtained four lists of respectively  711, 277, 319 and 369 probes correlated to the expression of ER$\alpha36$ in the TN, ERa66+,  PR+ and  PR- groups.

\noindent
{\bf Step 2b, functional analysis.} From these four lists of transcriptomic probes, we carried out a functional analysis using the MSigDB database (GSEA). In particular, we looked for transcription factors and microRNAs involved in the regulation of the majority of genes from the different lists (TN, ERa66+, PR+ and PR-).
The results indicated that four transcription factors: NFAT, FOXO4, SP1 and LEF-1, were common regulators and could therefore be mediators of the ER$\alpha36$ effect in all breast tumor subtypes. Interestingly, a study by    \cite{antanavivciute2017transcriptional} has shown that these transcription factors FOXO4, SP1 and LEF-1 are transcriptional hallmarks characteristic of cancer cells and associated to the Wnt signaling pathway (involved in metastasis and maintenance of cancer stem cells).
%Although these upstream factors were common to all tumor groups, when analyzing the Gene Ontology Biological Processes terms, we observed that these factors seemed to be involved in the regulation of different biological processes depending on the groups. Indeed, for TN and ER$\alpha66$+ tumors, the genes correlated with the expression of ER$\alpha36$ seemed mainly associated with metabolic processes whereas in the PR+ and PR- tumors we observed rather a link with the translational machinery and cell morphology, respectively. These terms remained very generic and it would be necessary to increase the number of transcriptomic data available to clarify the link between ER$\alpha36$ and these different biological processes.
Regarding the analysis of microRNAs, the results indicated that the majority of the genes whose expression correlated to ER$\alpha36$ one in the TN set were regulated by the microRNAs: hsa-miR-106B, hsa-miR21 and hsa-miR-29A, listed as oncogenic microRNAs involved in metastatic processes, survival and self-sufficiency in growth factors of mammary tumors \cite{le2010micro}. These results recalled those of a previous study of \cite{chamard2015eralpha66}, which showed that a high ER$\alpha36$ expression in mammary tumors is associated with an increased metastatic potential and an estrogen-independent tumor growth. %Conversely, in the PR- set, the genes correlated with ER$\alpha36$ expression were rather associated with micro-RNAs listed as tumor suppressors in mammary tumors: hsa-miR-34 and hsa- miR-200. These results indicate that there is a major post- transcriptional component in biological networks involving ER$\alpha36$.

\begin{figure}[htbp]
	\centering
		\includegraphics[scale=0.3]{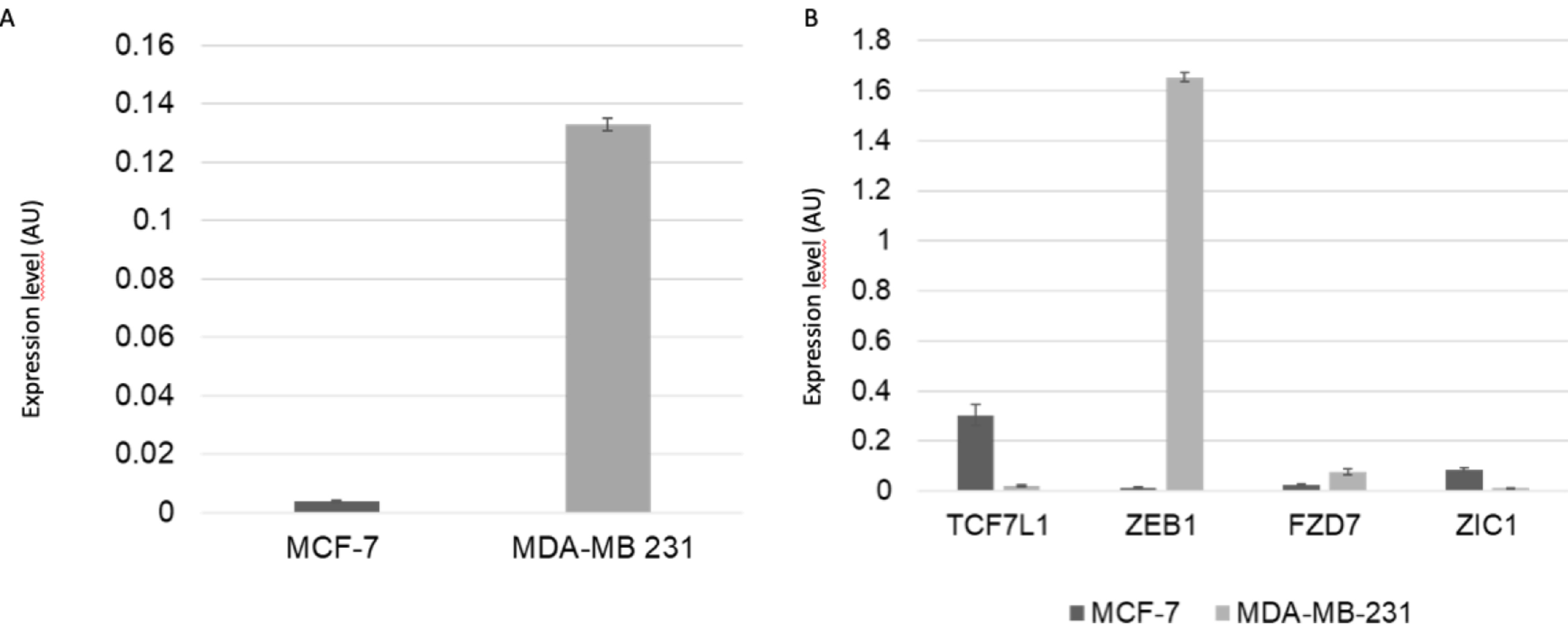}
		\caption{ Experimental validation of biological data inferred from \texttt{armada}.  (A) ER$\alpha36$ expression level as measured by RT-qPCR in MCF-7 and MDA-MB-231 breast cancer cells. (B) TCF7L1, ZEB1, FZD7 and ZIC1 expression level as measured by RT-qPCR in MCF-7 and MDA-MB-231 breast cancer cells.}
	\label{Fig2HD}
\end{figure}

\noindent
{\bf Step 2c, experimental validation.}
Subsequently, we provided an experimental confirmation of the biological reliability of the results: the correlation between the expression of ER$\alpha36$ and that of ZEB1, FZD7, ZIC1 and TCF7LD genes, identified by \texttt{armada} as correlated to ER$\alpha36$ in all tumor sets, was verified in vitro by RT-qPCR in two breast cancer cell lines (MCF-7 (ERa66+, PR+, PR-) and MDA-MB-231 (TN)). The results of Figure \ref{Fig2HD} confirmed the correlation (positive or negative) between the expression of ER$\alpha36$ and that of the genes identified by \texttt{armada} in the both cell lines.

\begin{figure}[htbp]
	\centering
		\includegraphics[scale=0.6]{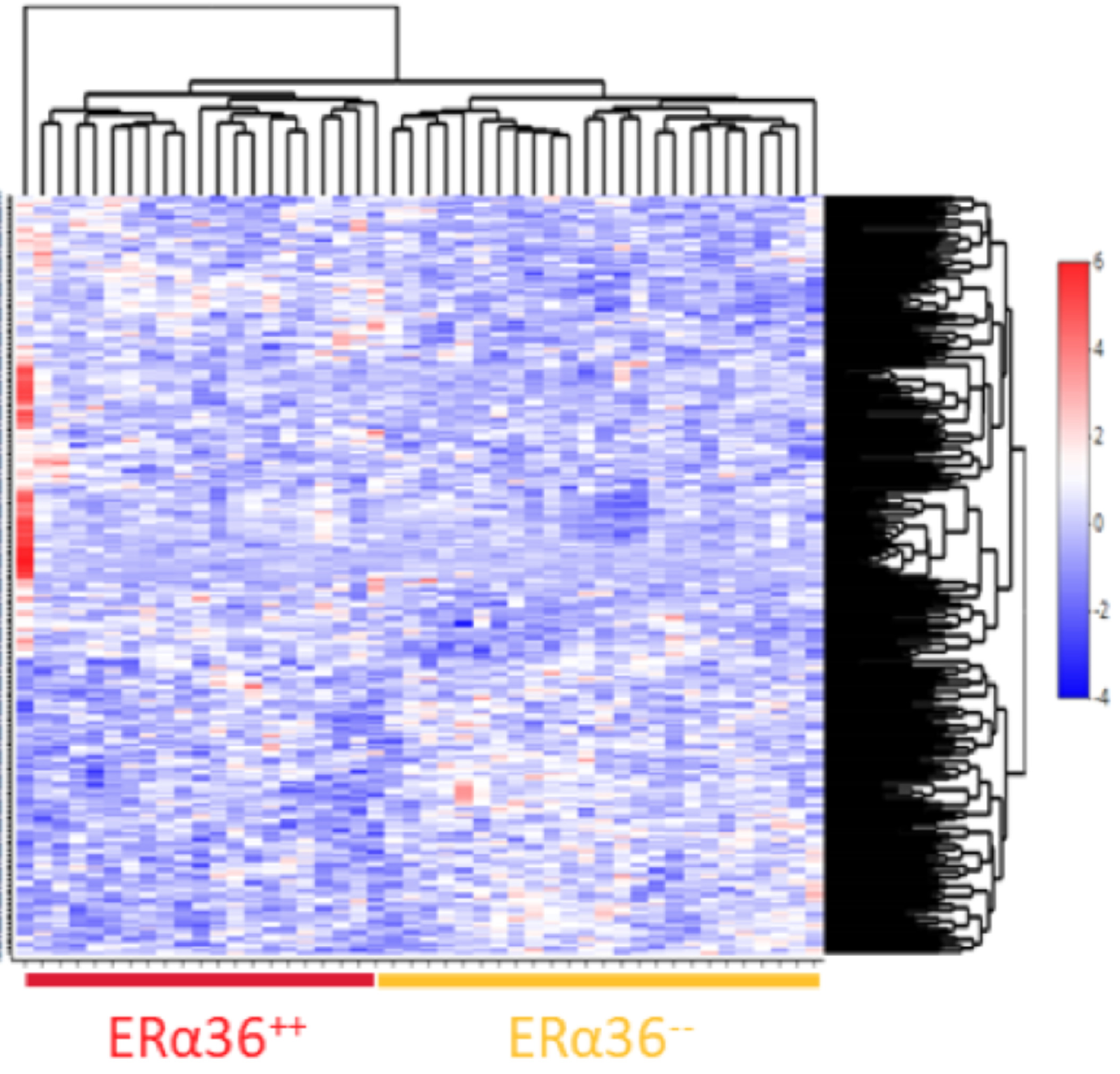}
		\caption{Heatmap of the 711 covariates which have ARMADA scores greater or equal to 1 in the study of the dataset ERa66+. Each column corresponds to one ERa66+ tumor. The $x$-axis represents the tumors (marked with their ER$\alpha36$ expression level) and the $y$-axis represents the selected probes. Color gives the normalised expression of the covariates.}
	\label{Fig3HD}
\end{figure}

\noindent
{\bf Step 3, tumor classification according to ER$\alpha36$ expression level.}
The final goal of our study was to identify the molecular signatures accounting for the ER$\alpha36$ expression level in the four different sets of tumors. These signatures were identified thanks to the R package \texttt{heatmaply} after co-clustering  of the ER$\alpha36$ expression (on the $x$-axis) and of the covariates (on the $y$-axis) with the function \texttt{hclust}. For each of the four tumor datasets, a heatmap was built which accounted for the expression level of the genes correlated to ER$\alpha36$. Thanks to the associated dendogram, different classes of tumors were defined and characterized by both the level of ER$\alpha36$ expression and an associated molecular signature. The Figure \ref{Fig3HD} illustrates the results for the study on the dataset ERa66+: two classes of tumors were identified, called ER$\alpha36^{++}$ and ER$\alpha36^-$. Taken together, the \texttt{armada} package helped to cluster patients which breast tumor highly express ER$\alpha36$ and associated genes. These patients could be treated by Wnt signaling inhibitors or specific microARN modulators and therefore benefit such promising new personalized medicine.

\section{Conclusion and perspectives}
\label{sec:conclu}

We have proposed a new methodology which is able to select the covariates (here the genes) that are linked with a variable of interest (here the treatment of an outcome or a biological marker). The method is of particular interest in the high dimensional case and when the covariates are correlated. The algorithms corresponding to this method are available through the  R package \texttt{armada}. After this selection obtained with our method, it is then easy to visualise the selected genes (or probes) for all the patients, and to classify the genetic profiles of patients with respect to their treatment outcome or biological marker. In the study of the treatment by chemotherapy in the advanced non-small-cell lung cancer, we have identify three types of genetic profiles defined with two clusters of genes.  In the study of the mammary tumors, the covariates selection allows the biologist to study the functional role of selected probes and also to classify tumors and associated transcriptomic signatures. %\begin{itemize}
%\item the patients who have an high level of expression for the genes of group G$_1$ and a low level of expression for the genes of group G$_2$ have all a low survival time (lower than 1 year);
%\item on the contrary, the patients who have a low level of expression for the genes of group G$_1$ and an high level of expression for the genes of group G$_2$ have all a better survival time (greater than 1 year);
%\item between these two opposite groups, we have a transitional group of patients who have intermediate expressions according to the two others groups. Apart from two patients, the survival time in this group is greater than 1 year (among the two exceptions, one patients has a survival time of 11.5 months).
%\end{itemize}
This kind of results is very promising for the identification of new therapeutic targets and the development of more efficient and personalized anti-cancer treatment.

\section*{Appendix}
The Appendix gives additional simulations on another classification model, and on a regression model where $Y$ is a continuous quantitative variable (Sections A and B of Appendix), some further analysis of the lung cancer data presented in Section \ref{sec:data1} (Section C of Appendix), and technical informations on the biological material used in Section \ref{sec:data2} (Section D of Appendix). 
\appendix

\section{Classification design}\label{sec:classif_Appendix}

\subsection{Simulation design}
As in Section \ref{sec:simu}, we propose a simulation study with $p=1600$ covariates and sample size $n=60$. The response variable $Y$ is binary: $Y=1$ for $\frac{n}{2}$ subjects, and $Y=0$ for $\frac{n}{2}$ subjects.  The covariates $\mathbf X=(\mathbf X^{(k)})_{k=1,\ldots,4}$ are clustered into four independent clusters, each of them containing $p_k=400$ covariates. For this, before to model the dependence with the outcome $Y$, we generate for each cluster $k$,  a preliminary vector $\mathbf{ \tilde{X}}^{(k)}$ that is a gaussian $400$-vector, with mean $0$ and non-diagonal variance-covariance matrix  $\mathbf{\Sigma}^{(k)}$. The correlation between the covariates of $\mathbf{\tilde{X}}^{(k)}$ inside the cluster $k$ is designed by a factor analysis model, as in the Section \ref{sec:simu}. More precisions on the factor analysis model can be found in  \cite{Friguet2012}. 
Now, we create the dependence between $\mathbf X$  and $Y$ in perturbing some component of $\mathbf{ \tilde{X}}$. This simulation design is inspired from the toys-data of \cite{genuer2010variable}. The outcome $Y$ is linked with 240 influential covariates in $\mathbf X$, the others being noise covariates. The links between the influential covariates  and the response variable $Y$ have different intensities. More precisely, the $10$ first covariates of each cluster are the most strongly linked with  the response variable $Y$  and the  strength of the link is decreasing in the successive groups of $10$ influential covariates.

More precisely, let us define the simulation model by giving the conditional distribution of $X_i$ given the value $y$ of  $Y$: in each cluster $k=1,\ldots, 4$, and for $i=1,\ldots, 400$, $${X}^{(k)}_i = \tilde{X}^{(k)}_i+\delta_i^{(k)}(y)$$ where $\delta_i(y)$ is a random variable.

\begin{itemize}
\item The relevant covariates are the $m_1=60$ first covariates of each cluster. The distribution of the $\delta_i^{(k)}(y)$ leading to the links between the relevant covariates and $Y$ is given in Table \ref{tab:toys_data_model}.
\item The $m_0 = 340$ remaining covariates of each cluster are independent of $Y$: $\delta_i^{(k)}(y)=0$ whatever $y$ for $i=61,\ldots,400$. 
\end{itemize}
\begin{table}[htbp]
\tbl{Links between the relevant covariates and $Y$ in the classification design. The notation $\delta_i^{(k)}\sim 0.7\mathcal N(3y, 1)+0.3\mathcal N(0,1)$ means that, with probability 0.7, $\delta_i^{(k)}\sim \mathcal N(3y, 1)$, and with probability 0.3, $\delta_i^{(k)}\sim \mathcal N(0,1)$.}
{\begin{tabular}{cc}\toprule
 $i$& model for $\delta_i^k$\\ \midrule
for $i=1,\ldots, 10$ & $\delta_i^{(k)}\sim 0.7\mathcal N(3y, 1)+0.3\mathcal N(0,1)$\\
for $i=11,\ldots, 20$ & $\delta_i^{(k)}\sim 0.7\mathcal N(2y, 1)+0.3\mathcal N(0,1)$\\
for $i=21,\ldots, 30$ & $\delta_i^{(k)}\sim 0.7\mathcal N(y, 1)+0.3\mathcal N(0,1)$\\
for $i=31,\ldots, 40$ & $\delta_i^{(k)}\sim 0.3\mathcal N(3y, 1)+0.7\mathcal N(0,1)$\\
for $i=41,\ldots, 50$ & $\delta_i^{(k)}\sim 0.3\mathcal N(2y, 1)+0.7\mathcal N(0,1)$\\
for $i=51,\ldots, 60$ & $\delta_i^{(k)}\sim 0.3\mathcal N(y, 1)+0.7\mathcal N(0,1)$\\
\bottomrule
\end{tabular}}
\label{tab:toys_data_model}
\end{table}

We can remark that this design respects the covariance matrix given in Figure \ref{fig:correl}. This design differs a little bit from the model of Equation (\ref{frigu}), because $\delta_i^{(k)}(Y)$ is a random function of  $Y$. Note that in real data analysis, we don't know the model from which they are generated. It is why it is interesting to analyse the performance of our method on different kinds of simulated data.

\subsection{Interest of our data pretreatment}

\noindent
In order to emphasize the interest of our data pretreatment, we compare the results of a Wilcoxon test after three different  data pretreatments:

\begin{enumerate}
\item[Procedure 1:] nothing is done on the dataset $\mathbf X$.
\item[Procedure 2:]  the  covariates $\mathbf X$ are decorrelated, taking $Y$ into account, with the factor analysis procedure of  \cite{friguet2009factor,causeur2011factor}, implemented in the R package  \texttt{FAMT}. This gives a new dataset $\mathbf X^\dagger_Y$.
\item[Procedure 3:] the 4 clusters are estimated with the procedure of \cite{chavent2011clustofvar}, implemented in the R package  \texttt{ClustOfVar};  then the covariates are decorrelated in each cluster,  taking $Y$ into account, with the factor analysis procedure of \cite{friguet2009factor,causeur2011factor}, implemented in the R package  \texttt{FAMT}. This gives a new dataset ${\mathbf X^*}_Y$ obtained by the  concatenation of the decorrelated clusters.
\end{enumerate}

\textit{\textbf{Remark:} our data pretreatment is the Procedure 3. We have supposed that the number of clusters is known. If that is not the case, the user can choose its own number of clusters by  using the graphical tools of the \texttt{ClustOfVar} procedure (plots of the dendrogram).}

\begin{figure}[htbp]
	\centering
		\includegraphics[scale=0.4]{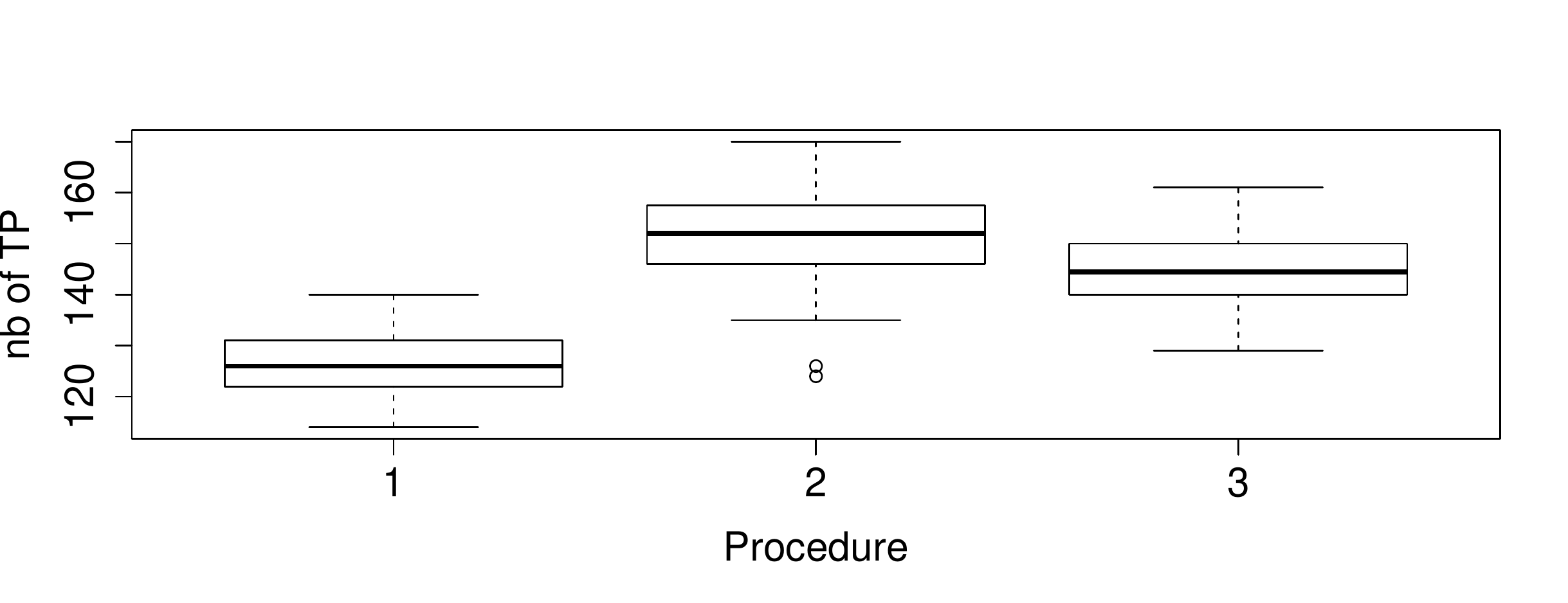}
		\includegraphics[scale=0.4]{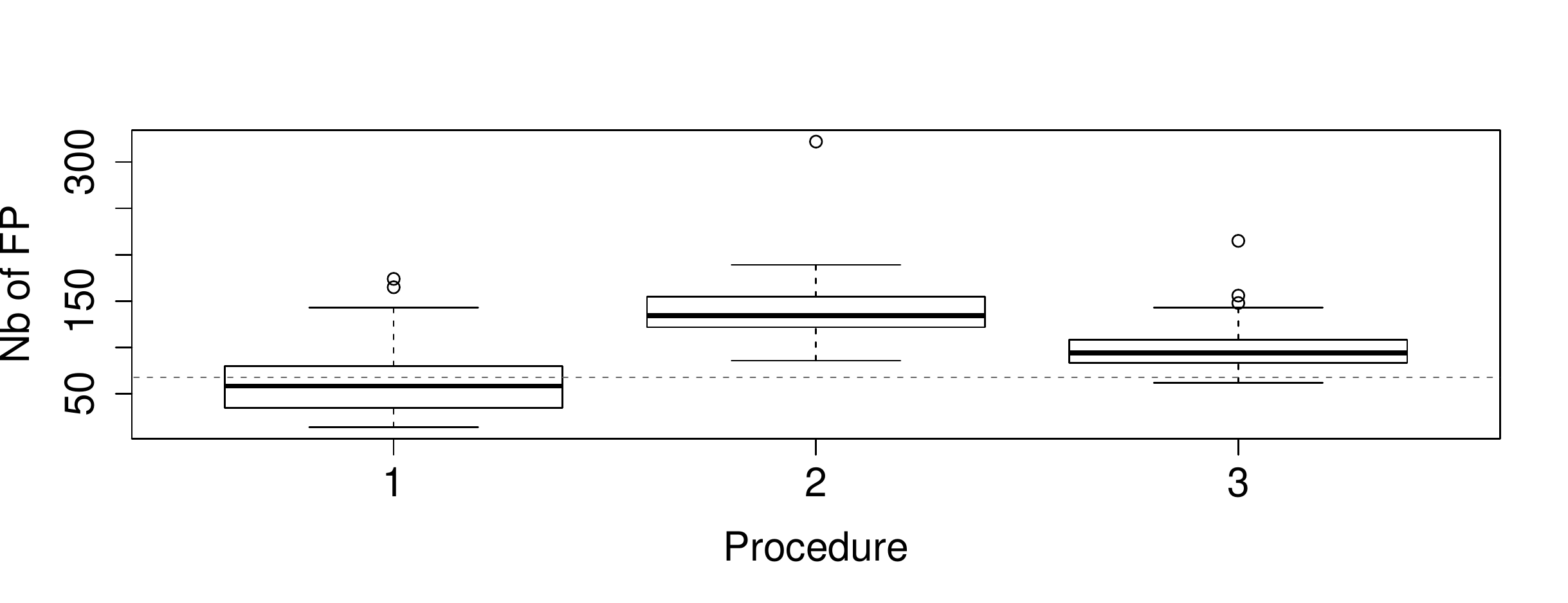}
				\caption{Number of true positive tests (top), false positive tests (bottom) in the classification design according to the different pretreatment procedures (1: Nothing, 2: FAMT, 3: clustering followed by FAMT in each cluster). Dotted lines: expected number of FP. Boxplots are calculated on $N=100$ runs.}
	\label{fig:nb_Positifs_delta2}
\end{figure}

Our objective is to find out the differently expressed covariates in the two groups (groups $Y=0$ and $Y=1$) with sample sizes $\frac{n}{2}=30$. For this, we perform Wilcoxon tests  on each of the $p$ pretreated covariates of the dataset (that is $\mathbf X$ for Procedure 1, $\mathbf X^\dagger_Y$ for Procedure 2, ${\mathbf X^*}_Y$ for Procedure 3), given a three sets of $p$ p-values. For each of these procedures, the selected covariates are those with p-values lower than  0.05. We compare these procedures on $N=100$ runs of $(\mathbf X,Y)$. For the comparison, we count the number of influential covariates that are correctly detected (this number is noted TP, for True Positive), this indicator gives an idea of the sensibility of the test after the procedure. To assess the specificity, we count the number of non-influential detected covariates (this number is noted FP, for False Positive). 
Note that the perfect method would detect all the influential covariates (that is 240 in this study) and no False Positive. However, according to the detection threshold chosen for the p-value, the expected number of FP is $68=5\%\times(1600-240)$. The results are shown in Figure \ref{fig:nb_Positifs_delta2}.

If we analyse the results given by Figure \ref{fig:nb_Positifs_delta2}, we can see that Procedure 1 is in fact the one that has the lowest rate of FP but its power is also the poorest whatever the design. Our Procedure reduces the mean and the variability of the distributions of the false positive rates, in comparison to the Procedure 2 (i.e. the FAMT procedure). The power of our Procedure is comparable with Procedure 2. This results show the interest of our proposed pretreatment before performing selection.

%%%%%%%%%%%%%%%%%%%%%%%%%%%%%
\subsection{Results of the whole method (pretreatment and selection)}
\begin{figure}[htbp]
	\centering
		\includegraphics[scale=0.4]{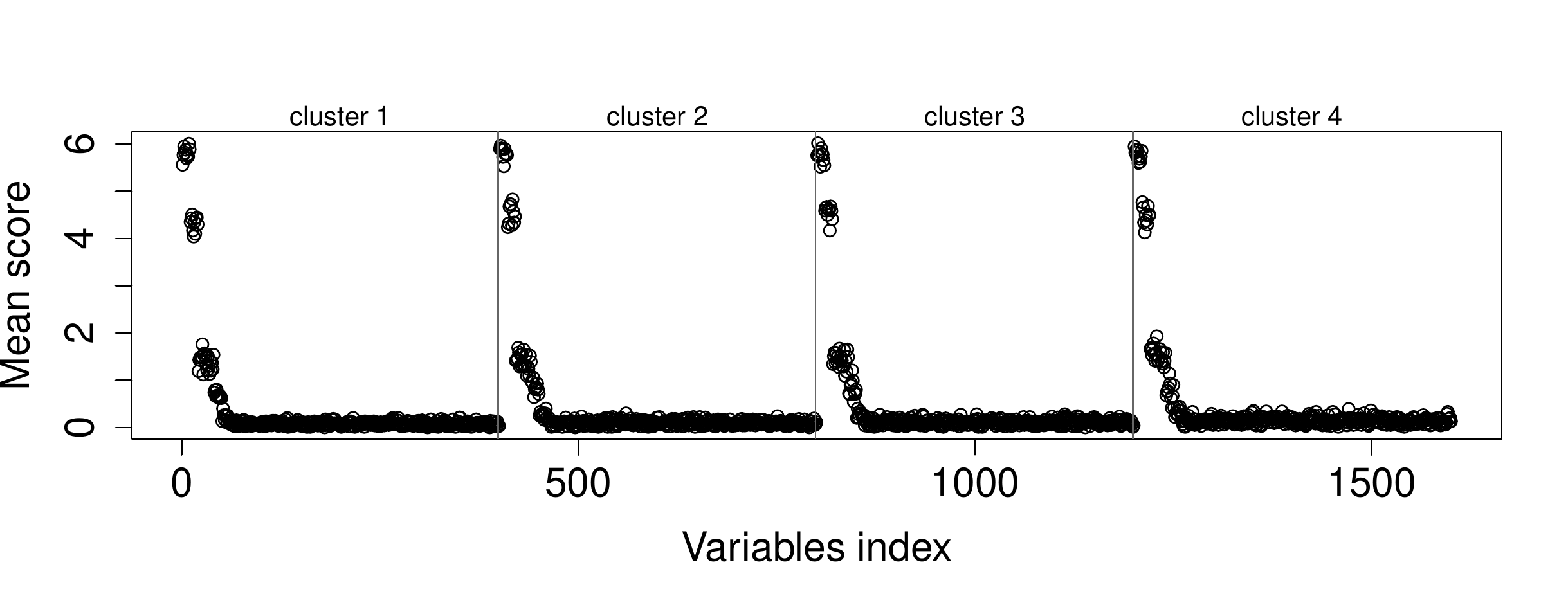}
	\includegraphics[scale=0.4]{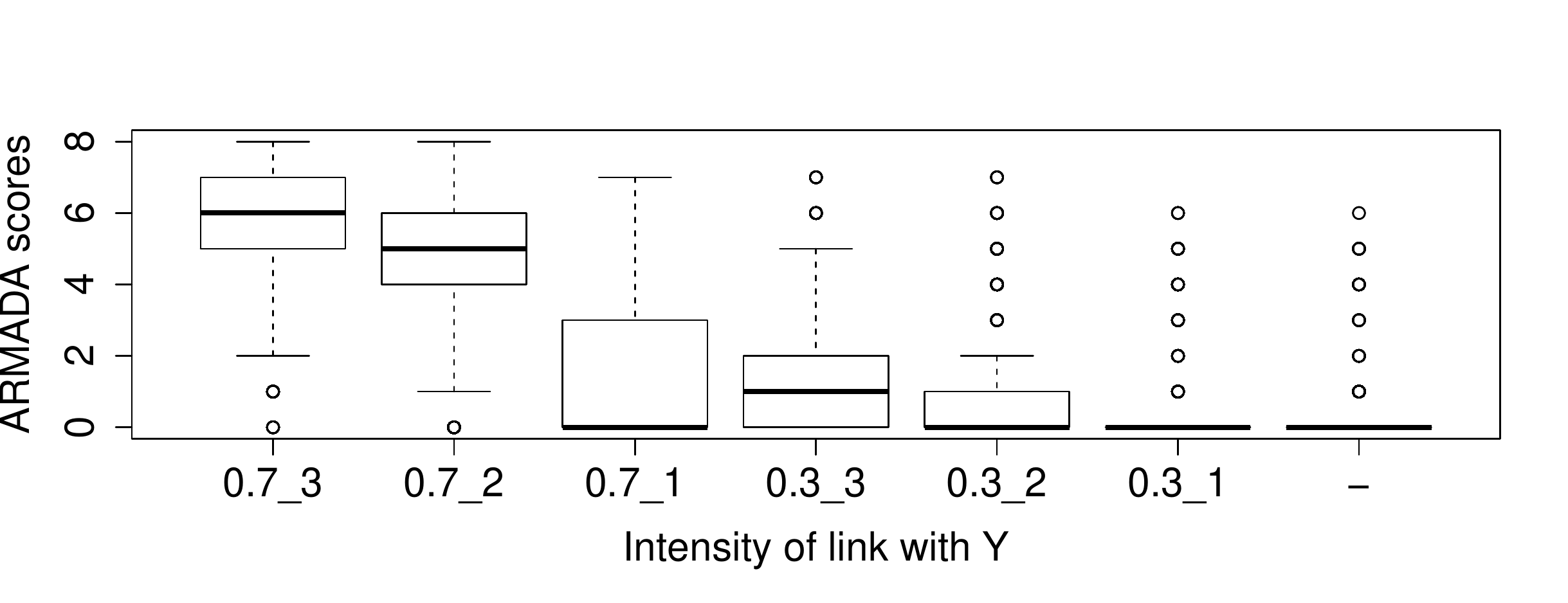}
			\caption{Top: mean of the ARMADA scores obtained by all the covariates. Bottom: boxplot of the scores of the covariates, ranked by levels of link with $Y$. Means and boxplots are calculated on $N=100$ runs. Simulation in the classification design.  }	\label{fig:mean_score_delta2}
\end{figure}

In order to describe the performances of our method, we show in Figure \ref{fig:mean_score_delta2}  the mean ARMADA scores obtained on the $N=100$ runs of $(\mathbf X, Y)$. The  scores are given for  all the covariates individually, and also by group of influential and noise covariates (the groups of influential covariates are noted by "(0.7,3)", "(0.7,2)", "(0.7,1)", etc.;  the group of noise covariates is noted by "-").

We can see on the Figure  \ref{fig:mean_score_delta2} that the scores give a clear ranking of the covariates, according to the strength of their link with the response variable $Y$. The highest scores are obtained by the covariates which are the most strongly linked with the response variable $Y$. 
The method is not so performant as in the design presented in Section \ref{sec:simu}, probably because we are note exactly in the model of the study, given by Equation (\ref{frigu}), but also because the strength of the link with $Y$ is low excepted for the two first groups of covariates that have scores which are well separated from the others by the selection method. 
We can precise that around 95\% of the noise covariates obtained an ARMADA score that was exactly 0.

\subsection{Comparison with other selection methods}\label{sec:comparison}
We propose the following selection criterion in our procedure: the selected covariates are those with scores greater or equal to 1. 

We compare this selection procedure with two other selection methods: 
\begin{itemize}
\item the Wilcoxon test: the selected covariates are those with raw-pvalues (i.e. p-values without any correction) lower than 0.05,
\item the FAMT procedure \cite{causeur2011factor}: the selected covariates are those with adjusted p-values lower than 0.05. 
\end{itemize}

To compare the three selection methods, the Table  \ref{tab:recap2_delta2} gives the rates of selection for each group of influential covariates, and for the group of noise covariates. The rates of selection have been computed on $N=100$ runs of $(\mathbf X,Y)$. We can see that our method respect the expected rate of false positives that is not the case for the FAMT method which exhibits a greater rate of 10 \%. 
 Our method is competitive with the FAMT procedure for the detection of influential covariates, but again FAMT procedure has more false positives than ours. 

\begin{table}[htbp]
\tbl{Results of the $N=100$ runs in the classification design: rates of selection of the different groups of influential and noise covariates by the ARMADA method, the Wilcoxon test and the FAMT procedure. The corresponding standard deviations are given in brackets. }
{\begin{tabular}{cccc}\toprule
 & ARMADA & Wilcoxon & FAMT\\ \midrule
(0.7-3)&                     0.99 (0.08) &       0.99 (0.07) &             0.99 (0.04)\\
(0.7-2)&                     0.92 (0.27) &                     0.92 (0.26) &                        0.96 (0.17)\\
(0.7-1) &                     0.44 (0.49)&                     0.43 (0.49) &                        0.58 (0.49)\\
(0.3-3)&                    0.54 (0.49)&                     0.41 (0.49) &                        0.61 (0.48)\\
(0.3-2)&                     0.32 (0.46) &                     0.28 (0.45) &                        0.41 (0.49)\\
(0.3-1)&                     0.12 (0.32)&                     0.12 (0.33)            &            0.19 (0.39)\\
-&                         0.05 (0.23)&                     0.05 (0.22) &                        0.09 (0.29)\\
\bottomrule
\end{tabular}}
\label{tab:recap2_delta2}
\end{table}

Finally, we can conclude with the ROC curves given in Figure \ref{fig:ROC}  that our method outperforms the two others selection methods (the ordinates of the points of the ARMADA ROC curve are all higher than the ordinates of the points of the two other ROC curves). Note that the ROC curves give the impression that our method is not competitive with the two others, but this is only caused by the fact that we have traced a solid line between the points $\textup{(1-specificity, sensibility)}_{\textup{ARMADA score} = 0}$ and $\textup{(1-specificity, sensibility)}_{\textup{ARMADA score} = 1}$. The ROC curves have been obtained by the mean of the $N=100$ ROC curves obtained in the $N=100$ runs of $(\mathbf X,Y)$.

 \begin{figure}[htbp]
	\centering
	\includegraphics[scale=0.4]{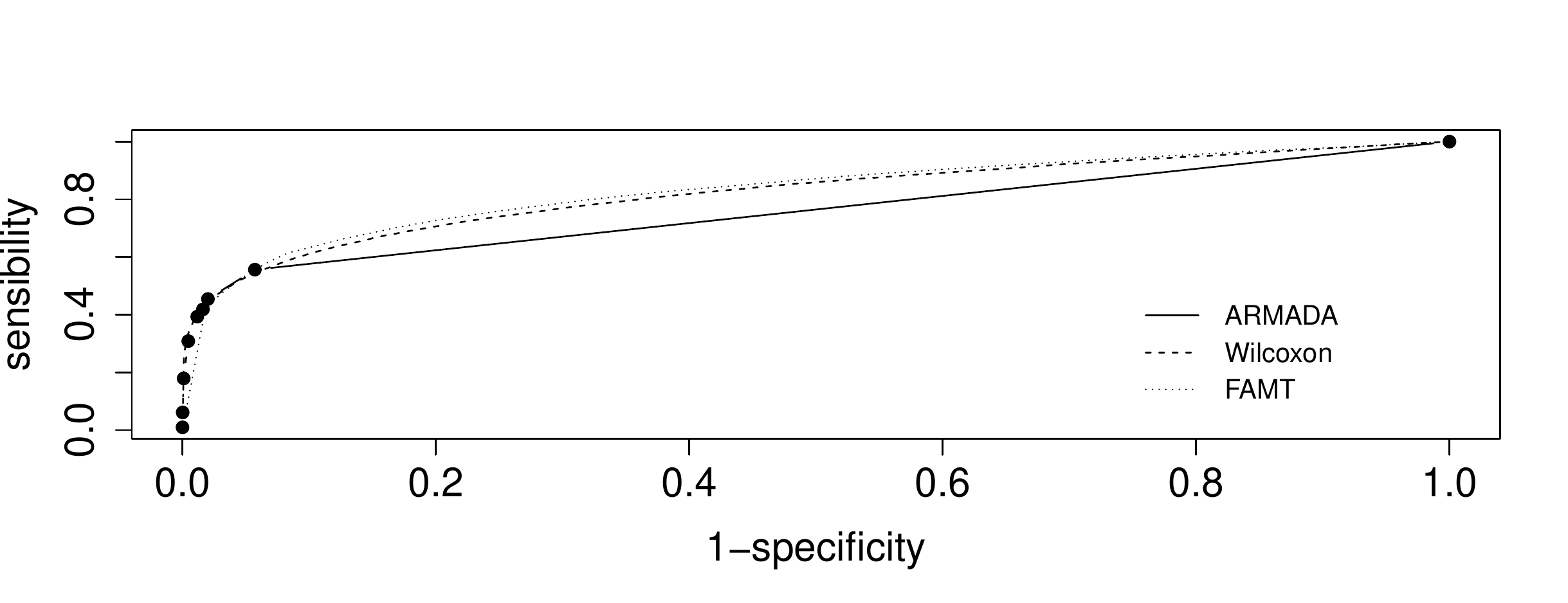}
			\caption{ROC curves for the three selection methods in the classification design. }	\label{fig:ROC}
\end{figure}

\section{Regression design}\label{sec:appendix}
In this section, we give results of simulations  to study the behavior of our algorithm to select covariates linked with a continuous variable of interest (like survival time here). 
We simulate $\mathbf{\tilde X}=(\mathbf{\tilde X}^{(k)})_{k=1,\ldots,4}$  as in Section \ref{sec:simu}, and $Y$ as a standard gaussian variable. Now, we create the dependence with outcome $Y$ in perturbing some component of $\mathbf{ \tilde{X}}$:
in all cluster $k=1,\ldots, 4$, and for all $i=1,\ldots, 400$: 
\begin{equation}
X_i^{(k)} = \tilde X_i^{(k)}+ \delta_iY
\end{equation}
where $\delta=\left(\delta_j\right)_{j=1,\ldots, 400}=(1,0.8,0.6,0.4,0.2, 0,0,\ldots, 0)$. Only the first 5 covariates of each cluster are linked with $Y$.

\begin{figure}[htbp]
	\centering
		\includegraphics[scale=0.4]{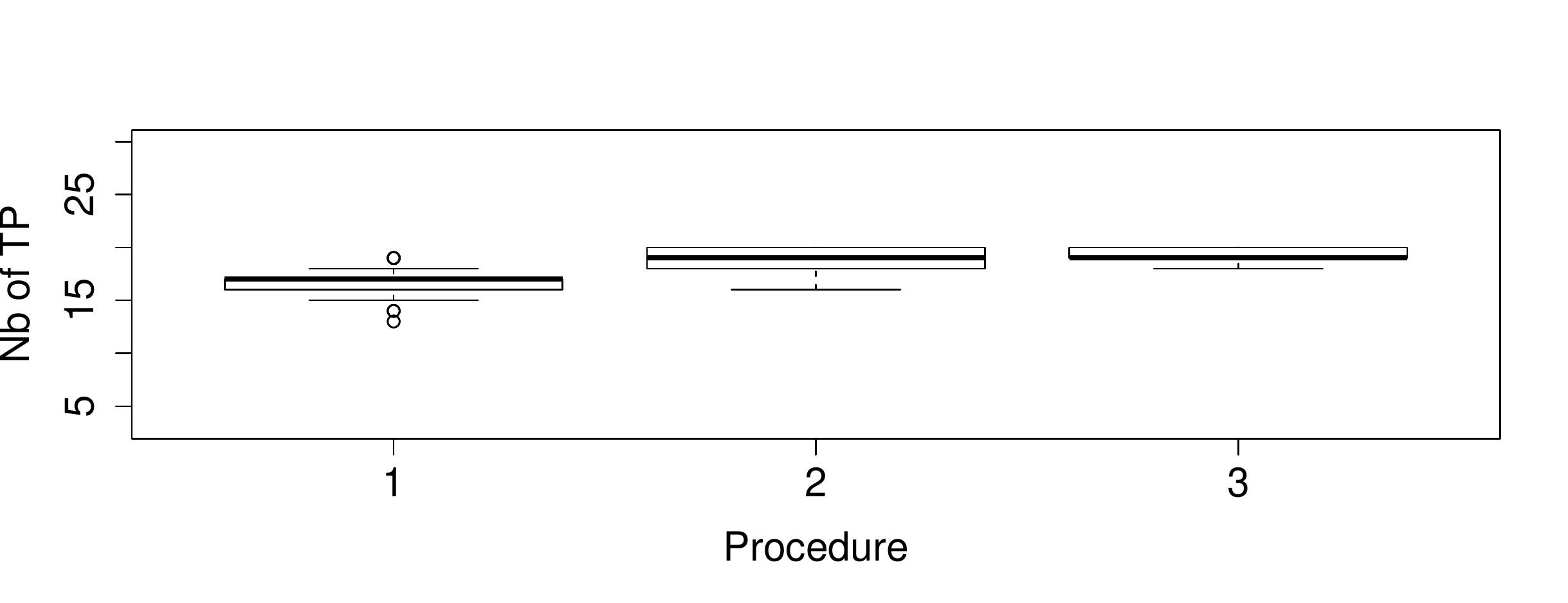}
		\includegraphics[scale=0.4]{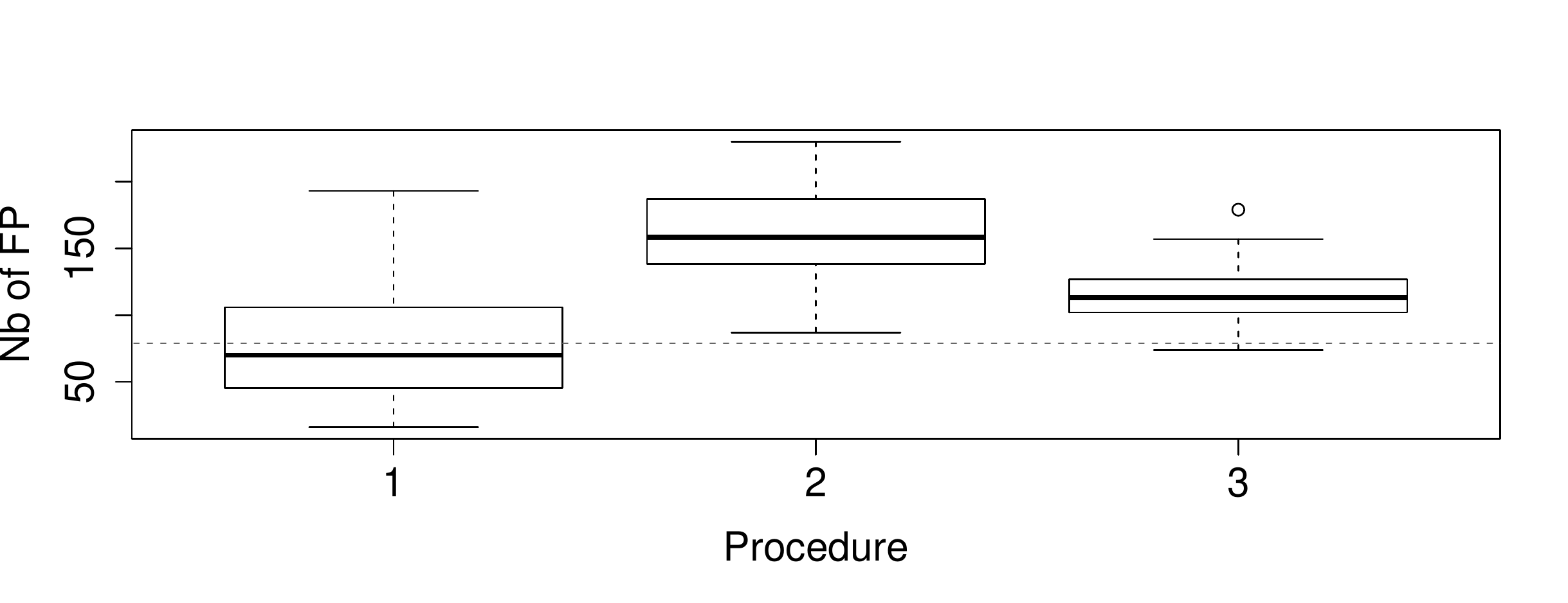}
				\caption{Number of:  true positive tests (top), false positive tests (bottom) in the regression design. Boxplots are calculated on $N=100$ runs.}
	\label{fig:nb_Positifs_reg}
\end{figure}

\begin{figure}[htbp]
	\centering
		\includegraphics[scale=0.4]{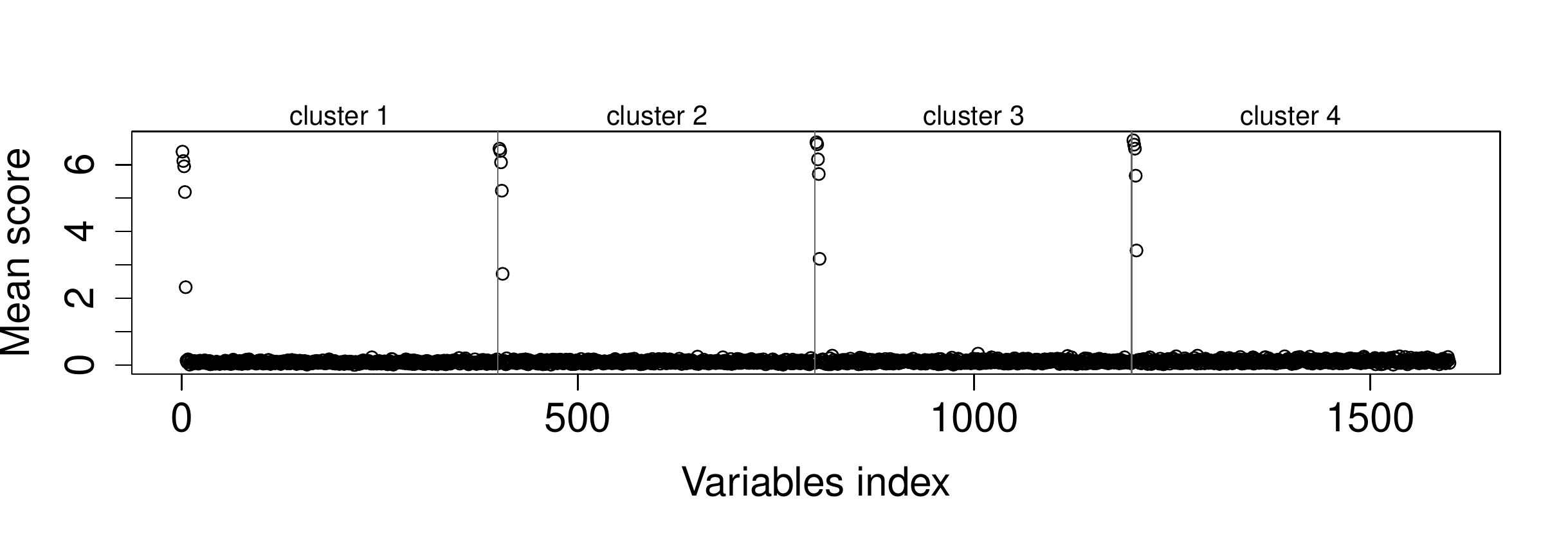}
	\includegraphics[scale=0.4]{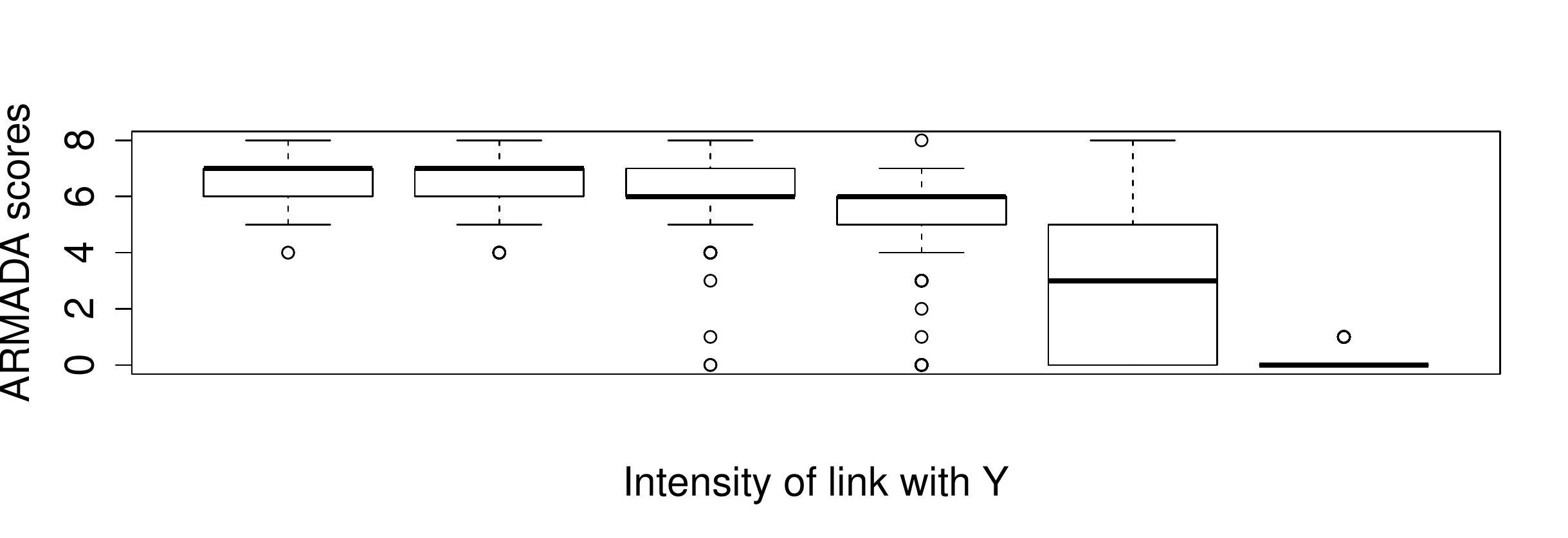}
			\caption{Top: mean of the ARMADA scores obtained by all the covariates. Bottom: boxplot of the scores of the covariates, ranked by levels of link with $Y$. Means and boxplots are calculated on $N=100$ runs. Simulation in the regression design. }	\label{fig:mean_score_reg}
\end{figure}

%%%%%%Pretraitement%%%%
We show the interest of our pretreatment, comparing the three procedures detailed in Section \ref{sec:simu}. As $Y$ is a gaussian variable, we use the Pearson correlation test (instead of the Wilcoxon test used in Section \ref{sec:simu}). We produce $N=100$ runs of $(\mathbf X, Y)$ and count the number of  false and true positive, and the ARMADA scores (shown in Figures \ref{fig:nb_Positifs_reg} and \ref{fig:mean_score_reg}). 

Similarly to the classification studies presented in Section \ref{sec:simu}, our Procedure reduces the mean and the variability of the distributions of the false positive rates, in comparison to the Procedure 2 (i.e. the FAMT procedure), and the power of our Procedure is comparable with Procedure 2.

The Figure \ref{fig:mean_score_reg} shows the ARMADA scores obtained on these $N=100$ runs of $(\mathbf X, Y)$.  Again, similarly to the Section \ref{sec:simu}, the scores give a ranking of the covariates, according to the intensity of their link with respect to the response variable $Y$. The true covariates are clearly separated of the noise covariates. We can also precise that 96\% of the noise covariates obtained a score that was 0.  
%R&\'esultas

\begin{figure}[htbp]
	\centering
		\includegraphics[scale=0.4]{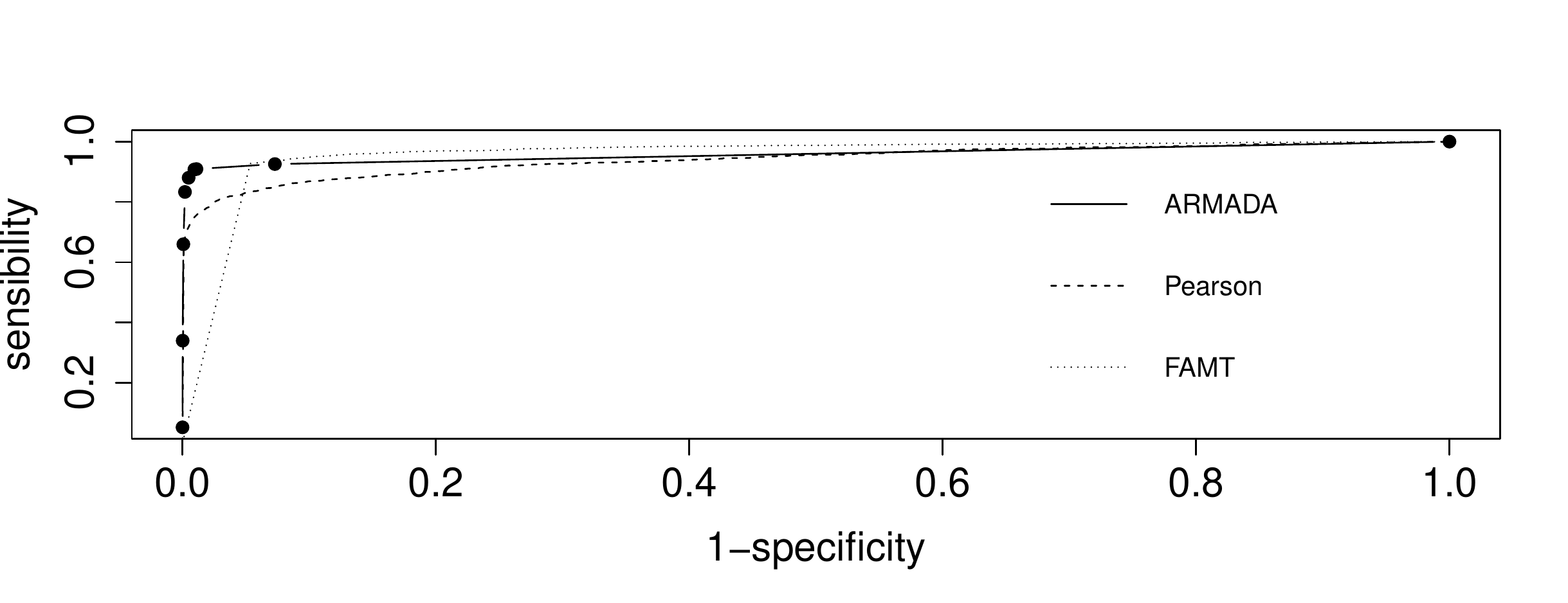}
			\caption{ROC curves for the three selection methods, in the case of regression design. The ROC curves have been obtained by the mean of the $N=100$ ROC curves obtained in the $N=100$ runs of $(\mathbf X,Y)$. }	\label{fig:ROCreg}
\end{figure}

%Comparison
As in Section \ref{sec:simu}, the Table \ref{tab:recap2_reg} and the ROC curve in Figure \ref{fig:ROCreg} allow us to compare our method with the Pearson test and the FAMT procedure. Our method seems to be a good compromise to have quite good detection rates for the true covariates, but small detection rates for the noise covariates. Even though true covariates are not always enough detected, compared to the FAMT procedure, detection rate of noisy covariates is lower than FAMT. The Pearson test has the lowest levels of detection rates, and the true covariates with a small link with $Y$ are not well detected. On the whole, our method seems to be appropriate for sparse models particularly when the goal is to avoid false positive detections.

\begin{table}[htbp]
\tbl{Results of the $N=100$ runs in the regression design: rates of selection of the different groups of influential and noise covariates by the ARMADA method, the Pearson correlation test and the FAMT procedure. The corresponding standard deviations are given in brackets.}
{\begin{tabular}{cccc}\toprule
 & ARMADA & Pearson & FAMT\\ \midrule
 1 &1 (0) &1 (0) &1 (0)\\
0.8&                     1 (0) &                     1 (0)&1 (0)\\
0.6&                     0.99 (0.08) &0.99 (0.08) &1 (0)\\
0.4&                     0.97 (0.18) &0.82 (0.38) &0.98 (0.13)\\
0.2&                     0.67 (0.47) &                     0.33 (0.47) &                       0.76 (0.43)\\
-&                     0.07 (0.26) &                     0.05 (0.22)&                       0.10 (0.30)\\
\bottomrule
\end{tabular}}
\label{tab:recap2_reg}
\end{table}

  \section{Lung cancer real dataset: bootstrap analysis}
\label{SM:transgene}
\begin{figure}[htbp]
	\centering
		\includegraphics[scale=0.5]{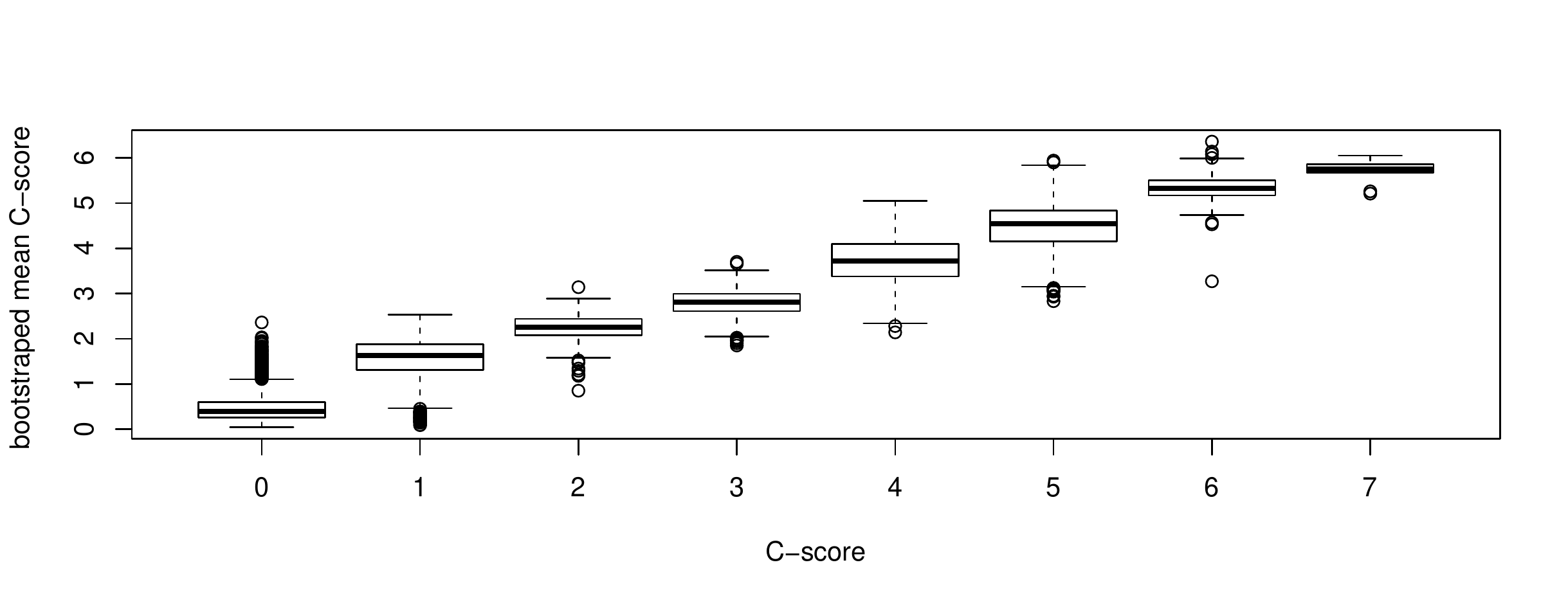}
		\includegraphics[scale=0.5]{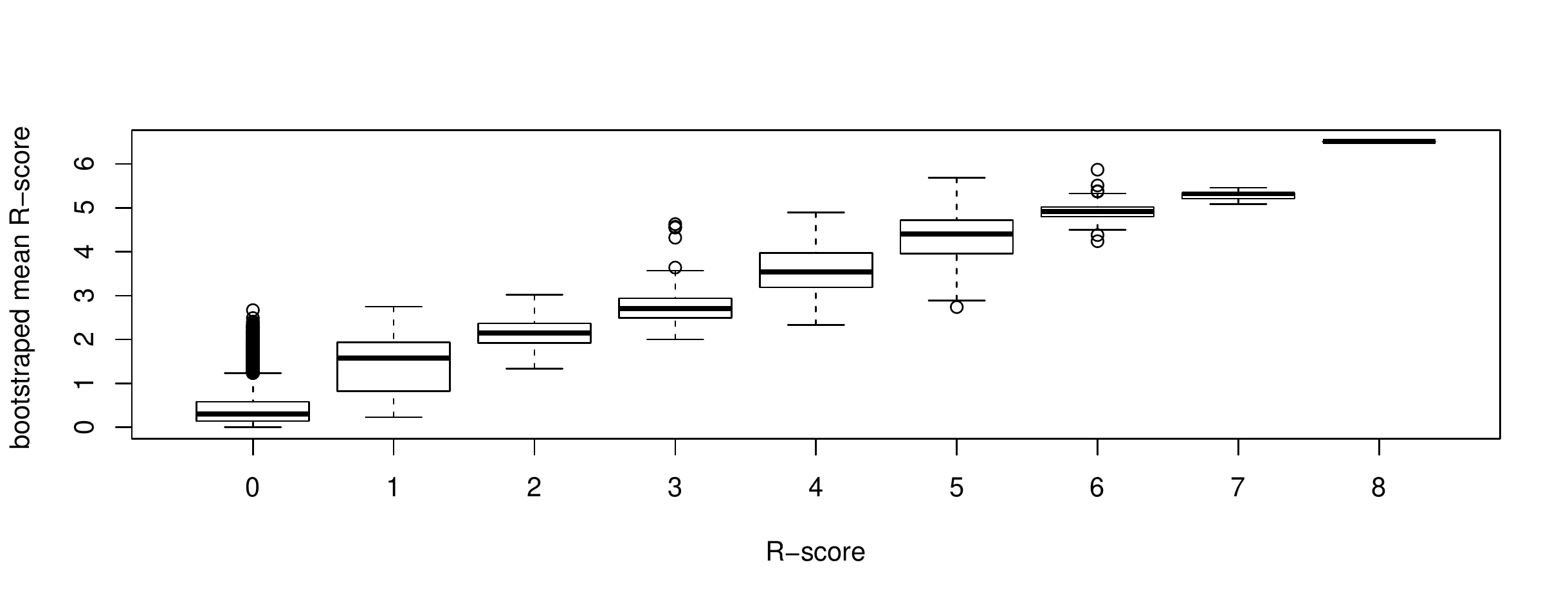}
	\caption{Distribution of the bootstraped mean of C-scores (resp. R-scores), i.e. means of C-(or R-)scores obtained on $B=100$ bootstrap samples), according to the corresponding C-scores (resp. R-scores) in the original dataset for all the $p=6810$ covariates.  }
	\label{fig:bootstrap_classif}
\end{figure}

As the number of patients $n=37$ is small compared to the number of covariates even after filtering ($p=6810$), we have checked our results with a bootstrap study. We have calculated the C-scores and R-scores of each covariates on $B=100$ bootstrap samples and the mean of the $B$ results.  We give the distribution of the bootstrapped means according to the original scores for the original dataset (Figure \ref{fig:bootstrap_classif}). 

We can see that the distributions of the bootstrapped means of the scores have a quite small dispersion and faithfully reproduce the original scores. The same conclusion holds for the bootstrapped median scores (shown in Tables \ref{tab:bootstrap_classif} and \ref{tab:bootstrap_reg}).

\begin{table}[htbp]
\tbl{Distribution of the boostraped median C-scores of the $p=6810$ covariates, obtained on $B=100$ boostrap samples, versus  the corresponding  C-scores.} 
{\begin{tabular}{c|ccccccccc}\toprule
 & \multicolumn{8}{c}{ARMADA C-score }\\ 
Bootstraped median C-score &      0&    1&    2&    3&    4&    5&    6&    7\\ \midrule
         0 &  2698&   53&    0&    0&    0&    0&    0&    0\\
         0.5&  11&    1&    0&    0&    0&    0&    0&    0\\
         1&    108&  315&   29&    1&    0&    0&    0&    0\\
         1.5&   0&   19&    5&    0&    0&    0&    0&    0\\
         2&      9&  162&  308&   76&    5&    0&    0&    0\\
         2.5&    0&    2&   28&   14&    0&    0&    0&    0\\
         3&      1&    1&   87&  321&   55&    2&    0&    0\\
         3.5&    0&    0&    1&   29&   12&    1&    0&    0\\
         4&      0&    0&    2&  155&  922&  218&    1&    0\\
         4.5&    0&    0&    0&    0&   19&    6&    0&    0\\
         5&      0&    0&    0&    0&  157&  644&  221&    2\\
         5.5&    0&    0&    0&    0&    0&    0&    6&    0\\
         6&      0&    0&    0&    0&    0&   17&   78  &8  \\ \bottomrule
\end{tabular}}
\label{tab:bootstrap_classif}
\end{table}

\begin{table}[htbp]
\tbl{Distribution of the boostraped median R-scores of the $p=6810$ covariates, obtained on $B=100$ boostrap samples, versus  the corresponding  R-scores.} 
{\begin{tabular}{c|cccccccccc}\toprule
 & \multicolumn{8}{c}{ARMADA C-score }\\
Bootstraped median R-score &      0&    1&    2&    3&    4&    5&    6&    7&8\\ \midrule
         0  & 3773   &29&   11&    0&    0&    0&    0&    0&    0\\
         0.5&   20&    2&    5&    0&    0&    0&    0&    0&    0\\
         1&     67&   22&   17&    0&    0&    0&    0&    0&    0\\
         1.5&    8&    1&    9&    0&    0&    0&    0&    0&    0\\
         2&    109&   32&  243&   40&    1&    0&    0&    0&    0\\
         2.5&    4&    0&   22&    8&    2&    0&    0&    0&    0\\
         3&      7&    3&  147&  295&   80&    2&    0&    0&    0\\
         3.5&    0&    0&    0&   14&   13&    2&    0&    0&    0\\
         4&      0&    0&    2&  149&  788&  210&    0&    0&    0\\
         4.5&    0&    0&    0&    0&   10&   14&    0&    0&    0\\
         5&      0&    0&    0&    3&   90&  462&   85&    5&    0\\
         5.5    &0&    0&    0&    0&    0&    1&    0&    0&    0\\
         6&      0&    0&    0&    0&    0&    1&    1&    0&    1\\
         \bottomrule
\end{tabular}}
\label{tab:bootstrap_reg}
\end{table}
Moreover, we can emphasis that our method is robust to detect the most important covariates (for instance, the 10 covariates that have a C-score equal to 7, or the 6 covariates that have an R-score greater than 7): their corresponding bootstraped means of scores are also high, and their corresponding bootstraped median scores are greater than  5.

\section{Biological material for the study of ER$\alpha36$ in breast cancer}
\label{SM:bio}

We analysed the biological network involving ER$\alpha36$ through the use of 4 sets of Affymetrix transcriptomic data obtained from breast tumors of different molecular subtypes: the triple negative (noted TN), ER66+, PR+ and PR- datasets:
\begin{itemize}
\item the TN dataset corresponding to Affymetrix transcriptomic comprehensive data from 17 patients derived xenografts (PDX) breast tumors was extracted from the Xentech$^{\text{TM}}$ database with the permission of Olivier D\'eas and Stefano Cairo (MTA CXT-295 Xentech SAS/University of Lorraine ; \cite{reyal2012molecular}). 
\item the 3 other datasets (46 tumors ER66+, 29 tumors PR+, 16 tumors PR-) were part of those from the Carte d'Identit\'e des Tumeurs Program (CIT) from the Ligue Nationale Contre le Cancer described in \cite{guedj2012refined}. Transcriptomic raw data were kindly provided by Aur\'elien De Reynies and Jacqueline M\'etral. One microgram of cDNAs from each tumor sample gathered at the Oncogenetics laboratory, INSERM U735, Institut Curie-H\^{o}pital-Centre Ren\'e Huguenin, St Cloud, France was also kindly provided by Ivan Bieche to measure ER$\alpha36$ expression.
\end{itemize}

The measurement of ER$\alpha$36 expression in each tumor (Step 1: clinical data completion) has been done as described in \cite{thiebaut2017mammary}. Total RNA extraction of PDX samples and qPCR analyses were performed. The following primers were used for qRT-PCR : GAPDH forward (Fw) 5'-TGC-ACC-ACC-AAC-TGC-TTA-GC -3', GAPDH reverse (Rev) 5'-GGC-ATG-GAC-TGT-GGT-CAT-GAG -3', ER$\alpha36$ forward (Fw) 5'- ATG-AAT-CTG-CAG-GGA-GAG-GA-3', ER$\alpha36$ reverse (Rev) 5'- GGC-TTT-AGA-CAC-GAG-GAA-ACC-3'. Assays were performed at least in triplicate, and the mean values were used to calculate expression levels, using the $\Delta \Delta$C(t) method referring to GAPDH housekeeping gene expression.

\section*{Notes on contributors}
The first case study of Section \ref{sec:data1} comes from Transgene team thanks to B. Bastien. T. Boukhobza, H. Dumond and C. Thiebaut conducted the  biological study of Section \ref{sec:data2} and performed the functional analysis that follows the covariates selection with \texttt{armada}. S. Cairo and O. D\'eas ;  XenTech, Genopole, 91000 Evry (France) provided the TN dataset; A. De Reynies, I. Bieche and J. M\'etral from the Carte d'Identit\'e des Tumeurs program provided the access to transcriptomic raw data and biological samples (ERa66+, PR+, PR- datasets). 
The statistical methodology has been developed by A. G\'egout-Petit and A. Muller-Gueudin, and thanks to Y. Shi and H. Chakir during their Master internships.

\end{document}